\documentclass[amsfonts,amstex,psfig]{elsart}
\usepackage{natbib}
\usepackage{amsmath}
\usepackage{amssymb}
\usepackage{graphicx}

\newtheorem{theorem}{Theorem}[section]

\newtheorem{proposition}[theorem]{Proposition}
\newtheorem{corollary}[theorem]{Corollary}
\newtheorem{algorithm2}[theorem]{Algorithm}
\newtheorem{remark}[theorem]{Remark}
\newtheorem{example}[theorem]{Example}

\newtheorem{definition}[theorem]{Definition}

\newenvironment{proof}[1]{
\trivlist \item[\hskip \labelsep{\bf #1}]}{\hfill\mbox{$\Box$} 
\endtrivlist}
\newcommand {\junk}[1]{}
\newcommand {\hide}[1]{}
\newcommand {\hideproof}[1]{}

\newcommand {\s}        {\mbox{\rm sign}}
\newcommand {\D}     {\mbox{\rm D}}

\newcommand {\Q}         {{\mathbb Q}}

\newcommand {\ZZ} {{\rm Z}}
\newcommand {\RR} {{\mathcal R}}
\newcommand {\re} {{\mathbb R}}
\newcommand {\la}   {{\langle}}
\newcommand {\ra}   {{\rangle}}
\newcommand {\eps} {{\varepsilon}}
\newcommand {\E} {{\rm Ext}}

\newcommand {\Ker}      {\mbox{\rm Ker}}

\newcommand{\R}{{\rm R}}
\newcommand {\Tot} {{\rm Tot}}
\newcommand {\level} {{\rm level}}
\newcommand {\ancestor} {{\rm ancestors}}
\newcommand {\F} {\mathbb Q}
\newcommand {\M} {{\mathcal M}}
\newcommand {\N} {{\mathcal N}}
\newcommand {\HH} {{\rm H}}
\newcommand {\LL} {{\rm L}}
\newcommand {\Ch} {{\rm C}}
\newcommand {\A} {{\mathbb A}}

\input{diagrams}
\newarrow{Line}{.}{.}{.}{.}{.}
\begin{document}

\begin{frontmatter}

\title{Computing the First Few Betti Numbers of Semi-algebraic Sets
in Single Exponential Time}

\author{Saugata Basu \thanksref{label2}}
\ead{saugata.basu@math.gatech.edu}
\thanks[label2]{Supported in part by NSF Career Award 0133597 and 
an Alfred P. Sloan Foundation Fellowship.}

\address{School of Mathematics,
Georgia Institute of Technology, Atlanta, GA 30332, U.S.A.}

\begin{abstract}
In this paper we describe an algorithm that takes as input a description
of a semi-algebraic set  $S \subset \R^k$,
defined by a Boolean formula with atoms of the form
$P > 0, P < 0, P=0$ for $P \in {\mathcal P} \subset \R[X_1,\ldots,X_k],$
and outputs the first $\ell+1$ Betti numbers of $S$,
$b_0(S),\ldots,b_\ell(S).$ The complexity of the algorithm is 
$(sd)^{k^{O(\ell)}},$ where
where $s  = \#({\mathcal P})$ and $d = \max_{P\in {\mathcal P}}{\rm deg}(P),$
which is singly  exponential in $k$ for $\ell$ any fixed constant. 
Previously, singly exponential time algorithms were known only for
computing the Euler-Poincar\'e characteristic, the zero-th and the first
Betti numbers.
\end{abstract}


\end{frontmatter}

\section{Introduction}
Let $\R$ be a real closed field and $S \subset \R^k$ a semi-algebraic set
defined by a Boolean formula with atoms of the form
$P > 0, P < 0, P=0$ for $P \in {\mathcal P} \subset \R[X_1,\ldots,X_k]$
(we call such a set a ${\mathcal P}$-semi-algebraic  set and the corresponding
formula a ${\mathcal P}$-formula).
It is well known \cite{O,OP,Milnor,Thom,B99,GV} 
that the topological complexity of $S$ 
(measured by the various Betti numbers of $S$) is bounded by $(sd)^{O(k)}$,
where $s  = \#({\mathcal P})$ and $d = \max_{P\in {\mathcal P}}{\rm deg}(P).$
Note that these bounds are singly exponential in $k$.
More precise bounds on the individual 
Betti numbers of $S$ appear in \cite{B03}.
Even though the Betti numbers of $S$ are bounded singly exponentially
in $k$, there is no known algorithm for producing a singly exponential
sized triangulation of $S$ (which would immediately imply a singly exponential
algorithm for computing the Betti numbers of $S$). In fact, 
designing a
singly exponential time algorithm for computing the Betti numbers of 
semi-algebraic sets is one of the outstanding open problems in algorithmic 
semi-algebraic geometry.
More recently, determining the exact complexity of computing the Betti numbers
of semi-algebraic sets has attracted the attention of computational 
complexity theorists \cite{BC}, who are interested in developing a theory
of counting complexity classes for the Blum-Shub-Smale model of real Turing
machines.

Doubly exponential algorithms (with complexity $(sd)^{2^{O(k)}}$)
for computing all the Betti numbers are known, 
since it is possible to obtain a triangulation of $S$ in doubly
exponential time using cylindrical algebraic decomposition 
\cite{Col,BPR03}.
In the absence of singly exponential time algorithms for computing
triangulations of semi-algebraic sets, algorithms with single exponential
complexity are known only for the problems of testing emptiness
\cite{R92,BPR95},
computing the zero-th Betti number 
(i.e. the number of semi-algebraically connected components of 
$S$) \cite{GV92,Canny93a,GR92,BPR99},
as well as the Euler-Poincar\'e characteristic of $S$ \cite{B99}.
Very recently a singly exponential time algorithm has been developed
for the problem of computing the first Betti number of a given semi-algebraic
set \cite{BPR04}.

In this paper we describe an algorithm,
which given a family
${\mathcal P} \subset \R[X_1,\ldots,X_k]$,
a ${\mathcal P}$-formula describing 
a ${\mathcal P}$-semi-algebraic set  $S \subset \R^k$,
and a number $\ell, 0 \leq \ell \leq k$ as input,
outputs the first $\ell$ Betti numbers of $S$.
For constant $\ell$, the complexity of the algorithm is singly exponential
in $k$.
We remark that using Alexander duality, 
we immediately get a singly exponential algorithm for 
computing the top $\ell$ Betti numbers too. However, the complexity
of our algorithm becomes doubly exponential if we want to compute the
middle Betti numbers of a semi-algebraic set using it.

There are two main ingredients in our algorithm for computing the
first $\ell$ Betti numbers of a given closed semi-algebraic set. The first
ingredient is a result proved in \cite{BPR04}, which enables us to
compute a singly exponential sized cover of the given semi-algebraic set
consisting of closed, contractible
semi-algebraic sets, in single exponential time.
The number and the degrees of the polynomials used to define the sets
in this cover are also bounded singly exponentially. 

The second ingredient, which is the main contribution of this paper, is 
an algorithm which uses the covering algorithm recursively and computes in
singly exponential time a complex whose cohomology groups are isomorphic to
the first $\ell$ cohomology groups of the input set. This complex is of
singly exponential size. 

The main result of the paper is the following.\\
\noindent
{\bf Main Result:} For any given $\ell$, there is an algorithm that
takes as input a ${\mathcal P}$-formula describing a semi-algebraic set 
$S \subset \R^k$,
and outputs $b_0(S),\ldots,b_\ell(S).$ 
The complexity of the algorithm is $(sd)^{k^{O(\ell)}}$, 
where $s  = \#({\mathcal P})$ and $d = \max_{P\in {\mathcal P}}{\rm deg}(P).$
Note that the complexity is singly exponential in $k$ for every fixed $\ell$.

The paper is organized as follows. In Section \ref{sec:math_prelim}, 
we recall some basic definitions from algebraic
topology and fix notations. In Section \ref{sec:covering}
we describe the construction of the complexes which allow us to compute the
the first $\ell$ Betti numbers of a given semi-algebraic set. In Section
\ref{sec:algo_prelim} we recall the inputs, outputs and complexities of a few
algorithms described in detail in \cite{BPR04}, which we use in our
algorithm. In Section \ref{sec:main} we describe our algorithm
for computing the first $\ell$ Betti numbers, prove its correctness as well
as the complexity bounds.
Finally in Section \ref{sec:practical} we comment on issues 
related to practical implementation.

\section{Mathematical Preliminaries}
\label{sec:math_prelim}

In this section, we recall some basic facts about semi-algebraic sets
as well as the definitions of
complexes and double complexes of vector spaces, and fix some
notations.

\subsection{Semi-algebraic sets and their cohomology groups}
Let $\R$ be  a real closed field. 
If ${\mathcal P}$ is a finite subset of
$\R [X_1, \ldots , X_k]$, we write the  set of zeros
of ${\mathcal P}$ in $\R^k$ as
$$
\ZZ({\mathcal P},\R^k)=\{x\in \R^k\mid\bigwedge_{P\in{\mathcal P}}P(x)= 0\}.
$$
We denote by 
$B(0,r)$ the open ball with center 0 and radius $r$.

Let ${\mathcal Q}$ and ${\mathcal P}$ be finite subsets of 
$\R[X_1,\ldots,X_k],$ $Z = \ZZ({\mathcal Q},\R^k),$ and
$Z_r = Z \cap B(0,r).$
A  {\em sign condition}  on
${\mathcal P}$ is an element of $\{0,1,- 1\}^{\mathcal P}$.
The {\em realization of the sign condition
$\sigma$ over $Z$}, $\RR(\sigma,Z)$, is the basic semi-algebraic set
\label{def:R(Z)}
$$
        \{x\in \R^k\;\mid\; \bigwedge_{Q \in {\mathcal Q}} Q(x)=0
     \wedge \bigwedge_{P\in{\mathcal P}} \s({P}(x))=\sigma(P) \}.
$$
The {\em realization of the sign condition
$\sigma$ over $Z_r$}, $\RR(\sigma,Z_r),$ is the basic semi-algebraic set
\label{def:R(Z_r)}
$
        \RR(\sigma,Z) \cap B(0,r).
$
For the rest of the paper, we fix an open
ball $B(0,r)$ with center 0 and radius $r$ big enough
so that, for  every sign condition $\sigma$,
$\RR(\sigma,Z)$ and $\RR(\sigma,Z_r)$ are homeomorphic.
This is always possible by the local conical structure at infinity
of semi-algebraic sets 
(\cite{BCR}, page 225).

A closed and bounded semi-algebraic
set $S \subset \R^k$ is semi-algebraically triangulable (see \cite{BPR03}),
and we denote by $\HH^i(S)$  the $i$-th simplicial cohomology group of $S$
with rational coefficients.  The groups  $\HH^i(S)$ are invariant under
semi-algebraic homeomorphisms and coincide with the corresponding
singular cohomology groups when $\R = \re$. We denote by $b_i(S)$ the
$i$-th Betti number of $S$ (that is, the dimension  of $\HH^i(S)$ as a vector
space), and $b(S)$ the sum $\sum_i b_i(S)$.
For a closed but not necessarily bounded semi-algebraic set $S \subset \R^k$,
we will denote by $\HH^i(S)$  the $i$-th simplicial cohomology group of 
$S \cap \overline{B(0,r)}$, where $r$ is sufficiently large. 
The sets $S \cap \overline{B(0,r)}$ are semi-algebraically homeomorphic
for all sufficiently large $r> 0$, by
the local conical structure at infinity of 
semi-algebraic sets, 
and hence this definition makes sense.

The definition of cohomology groups of arbitrary semi-algebraic sets in
$\R^k$ requires some care and several possibilities exist. In this paper,
we follow \cite{BPR03} and define the cohomology groups of 
realizations of sign conditions as follows. 

Let $\R$ denote a real closed field and 
$\R'$ a real closed field containing $\R$.
Given a semi-algebraic set
$S$ in ${\R}^k$, the {\em extension}
of $S$ to $\R'$, denoted $\E(S,\R'),$ is
the semi-algebraic subset of ${\R'}^k$ defined by the same
quantifier free formula that defines $S$.
The set $\E(S,\R')$ is well defined (i.e. it only depends on the set
$S$ and not on the quantifier free formula chosen to describe it).
This is an easy consequence of the transfer principle \cite{BPR03}.

Now, let $S \subset \R^k$ be a ${\mathcal P}$-semialgebraic set, where
${\mathcal P} = \{P_1,\ldots,P_s \}$ is a finite subset of $\R[X_1,\ldots,X_k].$ 
Let $\phi(X)$ be a quantifier-free formula defining $S$. 
Let $P_i = \sum_{\alpha} a_{i,\alpha}X^\alpha$ where the $a_{i,\alpha} \in \R.$
Let
$A = (\ldots,A_{i,\alpha},\ldots )$ denote the vector of variables 
corresponding 
to the coefficients of  the polynomials in the family ${\mathcal P},$
and let
$a = (\ldots,a_{i,\alpha},\ldots) \in \R^N$ denote the vector of
the actual coefficients of the polynomials in ${\mathcal P}$. 
Let $\psi(A,X)$ denote the formula obtained
from $\phi(X)$ by replacing each coefficient of each polynomial in ${\mathcal P}$
by the corresponding variable, so that $\phi(X) = \psi(a,X).$ It follows from
Hardt's triviality theorem for semi-algebraic mappings \cite{Hardt},
that there exists,
$a' \in \re_{\rm alg}^N$ such that 
denoting by $S' \subset \re_{\rm alg}^k$ the semi-algebraic set defined by
$\psi(a',X)$, the semi-algebraic set
$\E(S',\R)$ has the same homeomorphism type as $S$. 
Here, $\re_{\rm alg}$ is the field of real algebraic numbers.
 We define the 
cohomology groups of $S$ to be the singular cohomology groups 
of $\E(S',\re).$ It follows from the Tarski-Seidenberg transfer principle,
and the corresponding property of singular cohomology groups,
that the cohomology groups defined this way are invariant under 
semi-algebraic homotopies. It is also clear that this
definition is  compatible with the
simplicial cohomology for closed, bounded semi-algebraic sets, and 
the singular cohomology  groups when the ground field is $\re$.
Finally it is also clear that, the Betti numbers are not changed after 
extension:
\[
b_i(S)=b_i(\E(S,\R')).
\]

Note that we define the
co-homology groups of arbitrary semi-algebraic sets as above in order
to treat semi-algebraic sets over arbitrary (possibly non-archimedean)
real closed fields $\R$,
for which the standard proofs of the homology axioms 
(in particular the excision axiom) break down for singular homology groups
(see \cite{Knebusch}, page XIII).
If one is only interested in the case,
$\R = {\mathbb R}$, then singular co-homology groups suffice.

\subsection{Complex of Vector Spaces}
A sequence $\{\Ch^p\}$, $p \in {\mathbb Z}$, of
$\F$-vector spaces
together with a sequence
$\{\delta^p\}$ of
homomorphisms $\delta^p :\Ch^p \rightarrow \Ch^{p+1}$ 
(called differentials) for
which $\delta^{p-1}\; \delta^p = 0$ for all
$p$ is called a complex. 
When it is clear from context, we will drop the supercripts from the
differentials for the sake of readability.

The cohomology groups, $\HH^p(\Ch^\bullet)$ are defined by,
\[ 
\HH^p(\Ch^\bullet) = {Z^p(\Ch^{\bullet})}/{B^p(\Ch^{\bullet})},
\]
where
$B^p(\Ch^\bullet)= {\rm Im}(\delta^{p-1}),$
and 
$Z^p(\Ch^\bullet) = \Ker(\delta^p)$ and we will denote by
$\HH^*(\Ch^{\bullet})$ the graded vector space
$\bigoplus_{p} \HH^p(\Ch^{\bullet}).$

The cohomology groups, $\HH^p(\Ch^\bullet),$ are all
$\F$-vector spaces
(finite dimensional if the vector spaces $\Ch^p$'s are themselves finite
dimensional). We will henceforth omit reference to 
the field of coefficients $\F$ which is fixed throughout the rest of the
paper.

\subsection{Homomorphisms of Complexes}
Given two  complexes, $\Ch^\bullet = (\Ch^p,\delta^p)$ and $\D^{\bullet}=
(\D^{p},\delta^p)$,
a homomorphism of complexes,
$\phi^{\bullet}: \Ch^{\bullet} \rightarrow \D^{\bullet},$ is a
sequence of homomorphisms $\phi^p: \Ch^p \rightarrow \D^p$ for which
$\delta^p\; \phi^p = \phi^{p+1}\;\delta^p$ for all $p.$

In other words, the following diagram
is commutative.
\[
\begin{array}{ccccccc}
\cdots & \longrightarrow & \Ch^p &
\stackrel{\delta^p}{\longrightarrow} & \Ch^{p+1} & \longrightarrow &\cdots
\\ & &
\Big\downarrow\vcenter{\rlap{$\phi^p$}} & &
\Big\downarrow\vcenter{\rlap{$\phi^{p+1}$}} & & \\ \cdots & \longrightarrow
& \D^p &
\stackrel{\delta^p}{\longrightarrow} & \D^{p+1} & \longrightarrow &\cdots
\end{array}
\]

A homomorphism of complexes,
$\phi^{\bullet}: \Ch^{\bullet} \rightarrow \D^{\bullet},$ induces homorphisms,
$\phi^i: \HH^i(\Ch^{\bullet}) \rightarrow \HH^i(\D^{\bullet})$ and we will
denote the corresponding homomorphism between the  graded
vector spaces $\HH^*(\Ch^{\bullet}), \HH^*(\D^{\bullet})$ by
$\phi^*$. 
The homomorphism $\phi^{\bullet}$ is called a {\em quasi-isomorphism} if the 
homomorphism $\phi^*$ is an isomorphism.

Given two complexes $\Ch^{\bullet}$ and $\D^{\bullet}$, their direct sum
denoted by $\Ch^{\bullet} \oplus \D^{\bullet}$, is again a complex with its
$p$-th term being $\Ch^p \oplus \D^p$. Moreover, given two homomorphisms of
complexes,
$$
\displaylines{
\phi^{\bullet}: \Ch^{\bullet} \rightarrow \bar{\Ch}^{\bullet},\cr
\psi^{\bullet}: \D^{\bullet} \rightarrow \bar{\D}^{\bullet},
}
$$
their direct sum 
\[
\phi^{\bullet}\oplus\psi^{\bullet}:
\Ch^{\bullet} \oplus \D^{\bullet} \rightarrow
\bar{\Ch}^{\bullet} \oplus \bar{\D}^{\bullet},
\]
is again a homomorphism of complexes defined componentwise.
Note that if we specify a basis for the different terms of 
the complexes 
$\Ch^{\bullet},{\bar{\Ch}}^{\bullet},\D^{\bullet},\bar{\D}^{\bullet},$
as well as the matrices for the homomorphisms $\phi^{\bullet},\psi^{\bullet}$
then we can write down the matrix for the direct sum homomorphism
$\phi^{\bullet}\oplus\psi^{\bullet}$ as a sum of block-matrices using 
elementary linear algebra.

\subsection{The Nerve Lemma and Generalizations}
We first define formally the notion of a cover of a closed, bounded
semi-algebraic set.
\begin{definition}
\label{def:covering}
Let $S \subset \R^k$ be a closed and bounded semi-algebraic set.
A cover, ${\mathcal C}(S)$, of $S$ consists  of an ordered
index set, which by a slight abuse of language we also 
denote by  ${\mathcal C}(S)$, and a map that
associates to each $\alpha \in {\mathcal C}(S)$, 
a closed and bounded semi-algebraic
subset $S_\alpha \subset S$, such that 
$S = \cup_{\alpha \in {\mathcal C}(S)} S_\alpha$.
\end{definition}
For $\alpha_0,\ldots,\alpha_p, \in {\mathcal C}(S)$,
we associate to the formal product,
$\alpha_0\cdots\alpha_p$, the closed and bounded semi-algebraic set 
$S_{\alpha_0\cdots\alpha_p} = S_{\alpha_0} \cap \cdots \cap S_{\alpha_p}$.

Recall that the $0$-th simplicial cohomology group of a closed and bounded 
semi-algebraic set $X$, $\HH^0(X)$, can be identified with the 
$\Q$-vector space
of $\Q$-valued locally constant functions on $X$. Clearly, the dimension
of $\HH^0(X)$ is equal to the number of connected components of $X$.

For $\alpha_0,\alpha_1,\ldots,\alpha_p,\beta \in {\mathcal C}(S)$, 
and $\beta \not\in \{\alpha_0,\ldots,\alpha_p\}$,
let
\[
r_{\alpha_0,\ldots,\alpha_p;\beta}: 
\HH^0(S_{\alpha_0\cdots\alpha_p}) \longrightarrow 
\HH^0(S_{\alpha_0\cdots\alpha_{p}\cdot\beta})
\]
be the homomorphism defined as follows.
Given a locally constant function,
$\phi \in \HH^0(S_{\alpha_0\cdots\alpha_p})$,
$r_{\alpha_0\cdots\alpha_p;\beta}(\phi)$ 
is the locally constant function
on $S_{\alpha_0\cdots\alpha_{p}\cdot\beta}$ obtained by restricting
$\phi$ to $S_{\alpha_0\cdots\alpha_{p}\cdot\beta}$.

We define the generalized restriction homomorphisms,
\[
\delta^p: \bigoplus_{\alpha_0 < \cdots < \alpha_p, \alpha_i \in {\mathcal C}(S)}  
\HH^0(S_{\alpha_0\cdots\alpha_p})\longrightarrow 
\bigoplus_{\alpha_0< \cdots <\alpha_{p+1}, \alpha_i \in {\mathcal C}(S)} 
\HH^0(S_{\alpha_0\cdots\alpha_{p+1}})
\]
by 
\begin{equation}
\label{eqn:generalized_restriction}
\delta^p(\phi)_{\alpha_0\cdots\alpha_{p+1}}= \sum_{0 \leq i \leq p+1} 
(-1)^{i} r_{\alpha_0\cdots\hat{\alpha_i}\cdots\alpha_{p+1}; \alpha_i}
(\phi_{\alpha_0\cdots\hat{\alpha_{i}}\cdots\alpha_{p+1}}),
\end{equation}

where $\phi \in \bigoplus_{\alpha_0< \cdots <\alpha_p \in {\mathcal C}(S)}  
\HH^0(S_{\alpha_0\cdots\alpha_p})$ and
$r_{\alpha_0\cdots\hat{\alpha_i}\cdots\alpha_{p+1}; \alpha_i}$ 
is the restriction homomorphism defined previously.
The sequence of homomorphisms $\delta^p$ gives rise to a complex,
$\LL^{\bullet}({\mathcal C}(S))$, defined by,
\[
\LL^p({\mathcal C}(S)) = 
\bigoplus_{\alpha_0< \cdots <\alpha_p, \alpha_i \in {\mathcal C}(S)} 
\HH^0(S_{\alpha_0\cdots\alpha_p})
,
\]
with the differentials 
$\delta^p: \LL^p({\mathcal C}(S)) \rightarrow \LL^{p+1}({\mathcal C}(S))$ 
defined
in (\ref{eqn:generalized_restriction}).
The complex $\LL^{\bullet}({\mathcal C}(S))$ is 
often referred to as the {\em nerve complex} of the cover 
${\mathcal C}(S)$. 

For any $\ell \geq 0$,
we will denote by $\LL^{\bullet}_\ell({\mathcal C}(S))$ the truncated complex,
defined by,

$$
\begin{array}{cccc}
\LL^p_\ell({\mathcal C}(S)) &=& \LL^p({\mathcal C}(S)),& 0 \leq p \leq \ell, \\
&=& 0, & p > \ell.
\end{array}
$$

Notice that once we have a cover of $S$,
and we identify the connected components of the various intersections,
$S_{\alpha_0\cdots\alpha_p}$, we have natural bases for the vector
spaces
$$
\LL^p({\mathcal C}(S)) = 
\bigoplus_{\alpha_0< \cdots <\alpha_p, \alpha_i \in {\mathcal C}(S)} 
\HH^0(S_{\alpha_0\cdots\alpha_p})
$$ 
appearing as terms of the nerve complex. Moreover, the matrices
corresponding to the homomorphisms $\delta^p$ in this basis, depend only
on the inclusion relationships between the connected components of
$S_{\alpha_0\cdots\alpha_{p+1}}$ and those of 
$S_{\alpha_0\cdots\alpha_p}$. 

We say that the cover
${\mathcal C}(S)$
{\em satisfies the Leray property} if
each non-empty intersection
$S_{\alpha_0\cdots\alpha_p}$ is contractible.
Clearly, in this case
$$
\begin{array}{cccc}
\HH^0(S_{\alpha_0\cdots\alpha_p}) &\cong & \Q, &\mbox{if $S_{\alpha_0\cdots\alpha_p} \neq \emptyset$} \\
                                  &\cong &  0,  & \mbox{if $S_{\alpha_0\cdots\alpha_p} = \emptyset$}.
\end{array}
$$
It is a classical fact (usually referred to as the {\em nerve lemma}) that,

\begin{theorem}(Nerve Lemma)
\label{the:nerve}
Suppose that the cover ${\mathcal C}(S)$  satisfies the Leray property.
Then for each $i \geq 0$,
\[ \HH^i(\LL^{\bullet}({\mathcal C}(S))) \cong \HH^i(S).\]
\end{theorem}

\begin{proof}{Proof:}
See \cite{Rotman}.
\end{proof}

Thus, Theorem \ref{the:nerve} gives a method for computing the Betti numbers
of $S$ using linear algebra, 
from a cover of $S$ by contractible sets for which all
non-empty intersections are also contractible, once we are able to  test 
emptiness of the various intersections $S_{\alpha_0\cdots\alpha_p}$.

Now suppose that each individual member, $S_{\alpha_0}$ of the cover is 
contractible, but the various intersections  
$S_{\alpha_0\cdots\alpha_p}$
are not necessarily contractible for $p \geq 1$. 
Theorem \ref{the:nerve} does not hold in this case. 
However, the following is proved in \cite{BPR04}.

\begin{theorem}
\label{the:bettione}
Suppose that each individual member, $S_{\alpha_0}$ of the cover 
${\mathcal C}(S)$ is  contractible.
Then,
\[ \HH^i(\LL^{\bullet}_2({\mathcal C}(S))) \cong \HH^i(S),\]
for $i = 0,1.$
\end{theorem}

\begin{proof}{Proof:}
See \cite{BPR04}.
\end{proof}

Notice that from a cover by contractible sets, 
Theorem \ref{the:bettione} allows us to compute using linear algebra,
$b_0(S)$ and $b_1(S)$, once we have identified the non-empty 
connected components
of the pair-wise and triple-wise intersections of the sets in the cover,
and their inclusion relationships. It is quite easy to see that if we
extend the complex in Theorem  \ref{the:bettione} by one more term,
that is consider the complex, $\LL^{\bullet}_{3}({\mathcal C}(S))$,
then the cohomology of the complex does not yield information about
$\HH^2(S)$. Just consider the cover of the standard sphere 
$S^2 \subset \R^3,$  and the cover $\{H_1,H_2\}$ of $S^2$ where
$H_1,H_2$ are closed hemispheres meeting at the equator. The 
corresponding complex, $\LL^{\bullet}_{3}({\mathcal C})$,  is as follows.

\[
0 \rightarrow \HH^0(H_1) \bigoplus \HH^0(H_2) 
\stackrel{\delta^0} \longrightarrow 
\HH^0(H_{1}\cap H_{2}) \stackrel{\delta^1} \longrightarrow 
0
\longrightarrow
0
\]

Clearly, $\HH^2(\LL^{\bullet}_{3}({\mathcal C}(S))) \not\simeq \HH^2(S^2)$,
and indeed it is impossible to compute  $b_i(S)$ 
just from the information on the number of connected components of 
intersections of the sets of a cover by contractible sets for, $i \geq 2$.
For example, the nerve complex coresponding to the cover of the
sphere by two hemispheres is ismorphic to the nerve complex of a
cover of the unit segment $[0,1]$ by the subsets $[0,1/2]$ and 
$[1/2,1]$, but clearly $\HH^2(S^2) = \Q$, while $\HH^2([0,1]) = 0$. 

In order to deal with covers not satisfying the Leray property, it
is necessary to consider a generalization of the nerve complex, namely 
a double complex arising from the generalized Mayer-Vietoris exact sequence.
The construction of this double complex (which is quite classical)
in fact motivates the design of our algorithm, which we describe 
in detail in Section \ref{sec:main}. 

\section{Mayer-Vietoris}
\subsection{Double Complexes}
\label{sec:double}
In this section, we recall  the basic notions of a 
double complex of vector spaces and associated spectral 
sequences.
A {\em double complex} is a bi-graded vector space,
\[
{\Ch}^{\bullet,\bullet} = \bigoplus_{p,q \in {\mathbb Z}} \Ch^{p,q},
\]
with co-boundary operators
$d : \Ch^{p,q} \rightarrow \Ch^{p,q+1}$ and
$\delta: \Ch^{p,q} \rightarrow \Ch^{p+1,q}$ and such that $d\delta +\delta d = 0$.
We say that $\Ch^{\bullet,\bullet}$ is a first quadrant
double complex, if it satisfies
the condition that $\Ch^{p,q} = 0$ if either $p < 0$ or $q < 0$.
Double complexes lying in other quadrants are defined in an analogous manner.

\hide{
{\small
\begin{diagram}
\vdots &              & \vdots&               &\vdots&              & \\
\uTo^{d} &              & \uTo^{d}&               &\uTo^{d}&              & \\
\Ch^{0,2}  & \rTo^{\delta}& \Ch^{1,2} & \rTo^{\delta} & \Ch^{2,2}& \rTo^{\delta}&\cdots \\
\uTo^{d} &              & \uTo^{d}&               &\uTo^{d}&              & \\
\Ch^{0,1}  & \rTo^{\delta}& \Ch^{1,1} & \rTo^{\delta} & \Ch^{2,1}& \rTo^{\delta}&\cdots \\
\uTo^{d} &              & \uTo^{d}&               &\uTo^{d}&              & \\
\Ch^{0,0}  & \rTo^{\delta}& \Ch^{1,0} & \rTo^{\delta} & \Ch^{2,0}& \rTo^{\delta}&\cdots \\
\end{diagram}
}
}

The complex defined by 
\[
{\rm Tot}^n(\Ch^{\bullet,\bullet}) = \bigoplus_{p+q=n} \Ch^{p,q},
\]
with differential
\[\D^n  = d \pm \delta: {\rm Tot}^{n}(\Ch^{\bullet,\bullet}) 
\longrightarrow {\rm Tot}^{n+1}(\Ch^{\bullet,\bullet}),
\]
is denoted by ${\rm Tot}^{\bullet}(\Ch^{\bullet,\bullet})$ and called
the {\em associated total complex of $\Ch^{\bullet,\bullet}$}.

\hide{
{\small
\begin{diagram}
&\vdots &              & \vdots&               &\vdots&              & \\
\luLine&\uTo^{d} &  \luLine            & \uTo^{d}&  \luLine             &\uTo^{d}& \luLine              & \\
\rTo^{\delta}& \Ch^{p-1,q+1} &  \rTo^{\delta}&\Ch^{p,q+1}&\rTo^{\delta}&\Ch^{p+1,q+1}&\rTo^{\delta}&\cdots \\
\luLine&\uTo^{d} &  \luLine            & \uTo^{d}&  \luLine             &\uTo^{d}&        \luLine      & \\
\rTo^{\delta}&\Ch^{p-1,q}  & \rTo^{\delta}& \Ch^{p,q} & \rTo^{\delta} & \Ch^{p+1,q}& \rTo^{\delta}&\cdots \\
\luLine&\uTo^{d} &  \luLine            & \uTo^{d}&   \luLine            &\uTo^{d}&            \luLine  & \\
\rTo^{\delta}& \Ch^{p-1,q-1}  & \rTo^{\delta}& \Ch^{p,q-1} & \rTo^{\delta} & \Ch^{p+1,q-1}& \rTo^{\delta}&\cdots \\
\luLine&\uTo^{d} &   \luLine           & \uTo^{d}& \luLine               &\uTo^{d}&           \luLine   & \\
&\vdots &              & \vdots&               &\vdots&              & \\
\end{diagram}
}
}

\subsection{Spectral Sequences}
A {\em spectral sequence} is a sequence of bigraded complexes
$(E_r, d_r: E^{p,q}_r \rightarrow E^{p+r,q-r+1}_r)$ (see Figure 1)
such that 
the complex $E_{r+1}$ is obtained from $E_r$ by taking its cohomology
with respect to $d_r$ (that is $E_{r+1} = \HH_{d_r}(E_r)$).

{\small
\begin{figure}[hbt] 
\begin{center}
\begin{picture}(0,0)%
\includegraphics{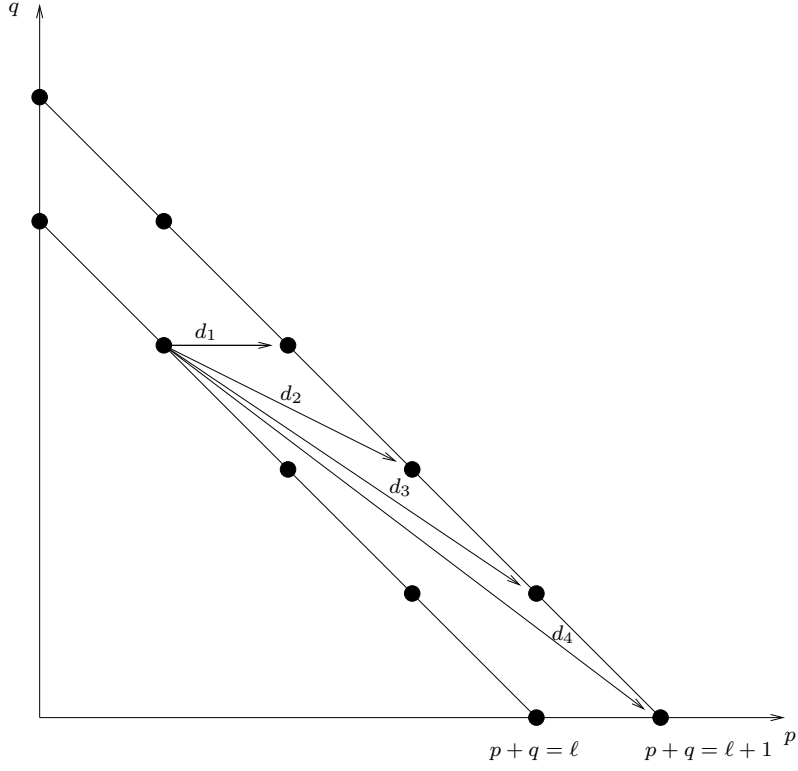}%
\end{picture}%
\setlength{\unitlength}{2565sp}%
\begingroup\makeatletter\ifx\SetFigFont\undefined%
\gdef\SetFigFont#1#2#3#4#5{%
  \reset@font\fontsize{#1}{#2pt}%
  \fontfamily{#3}\fontseries{#4}\fontshape{#5}%
  \selectfont}%
\fi\endgroup%
\begin{picture}(7512,7414)(601,-7394)
\put(6751,-7336){\makebox(0,0)[lb]{\smash{\SetFigFont{8}{9.6}{\familydefault}{\mddefault}{\updefault}{$p + q = \ell+1$}%
}}}
\put(5251,-7336){\makebox(0,0)[lb]{\smash{\SetFigFont{8}{9.6}{\familydefault}{\mddefault}{\updefault}{$p + q = \ell$}%
}}}
\put(8101,-7186){\makebox(0,0)[lb]{\smash{\SetFigFont{8}{9.6}{\familydefault}{\mddefault}{\updefault}{$p$}%
}}}
\put(601,-136){\makebox(0,0)[lb]{\smash{\SetFigFont{8}{9.6}{\familydefault}{\mddefault}{\updefault}{$q$}%
}}}
\put(2401,-3286){\makebox(0,0)[lb]{\smash{\SetFigFont{8}{9.6}{\familydefault}{\mddefault}{\updefault}{$d_1$}%
}}}
\put(3226,-3886){\makebox(0,0)[lb]{\smash{\SetFigFont{8}{9.6}{\familydefault}{\mddefault}{\updefault}{$d_2$}%
}}}
\put(4276,-4786){\makebox(0,0)[lb]{\smash{\SetFigFont{8}{9.6}{\familydefault}{\mddefault}{\updefault}{$d_3$}%
}}}
\put(5851,-6211){\makebox(0,0)[lb]{\smash{\SetFigFont{8}{9.6}{\familydefault}{\mddefault}{\updefault}{$d_4$}%
}}}
\end{picture}
\caption{$d_r: E_r^{p,q} \rightarrow E_r^{p+r, q- r +1}$}
\label{fig:spectral}
\end{center}
\end{figure}
}

There are two spectral sequences,
$'E_*^{p,q},{''E}_*^{p,q}$,  (corresponding to taking row-wise or
column-wise filtrations respectively) 
associated with a first quadrant 
double complex ${C}^{\bullet,\bullet},$ 
which  will be important for us. 
Both of these converge to $\HH^*({\rm Tot}^{\bullet}(\Ch^{\bullet,\bullet})).$
This means that the homomorphisms, $d_r$ are eventually zero, and hence the
spectral sequences stabilize, and
\[
\bigoplus_{p+q = i}{'E}_{\infty}^{p,q} \cong 
\bigoplus_{p+q = i}{''E}_{\infty}^{p,q} \cong 
\HH^i({\rm Tot}^{\bullet}(\Ch^{\bullet,\bullet})),
\]
for each $i \geq 0$.

The first terms of these are
\[{'E}_1 = \HH_{d}(\Ch^{\bullet,\bullet}), 
{'E}_2 = \HH_d\HH_{\delta}(\Ch^{\bullet,\bullet}),
\]
and
\[
{''E}_1 = \HH_{\delta} (\Ch^{\bullet,\bullet}), {''E}_2 = 
\HH_d \HH_{\delta} (\Ch^{\bullet,\bullet}).
\]

Given  two (first quadrant) double complexes, $\Ch^{\bullet,\bullet}$ and
$\bar{C}^{\bullet,\bullet},$ a homomorphism of double complexes,
\[
\phi^{\bullet,\bullet}: \Ch^{\bullet,\bullet} \longrightarrow \bar{C}^{\bullet,\bullet},
\]
is a collection of homomorphisms,
$\phi^{p,q}: \Ch^{p,q} \longrightarrow \bar{C}^{p,q},$ such that
the following diagrams commute.

\[
\begin{array}{ccc}
\Ch^{p,q} &
\stackrel{\delta}{\longrightarrow} & \Ch^{p+1,q} \\ 
\Big\downarrow\vcenter{\rlap{$\phi^{p,q}$}} & &
\Big\downarrow\vcenter{\rlap{$\phi^{p+1,q}$}} \\ 
\bar{C}^{p,q} &
\stackrel{\delta}{\longrightarrow} & \bar{C}^{p+1,q}
\end{array}
\]

\[
\begin{array}{ccc}
 \Ch^{p,q} &
\stackrel{d}{\longrightarrow} & \Ch^{p,q+1}\\
\Big\downarrow\vcenter{\rlap{$\phi^{p,q}$}} & &
\Big\downarrow\vcenter{\rlap{$\phi^{p,q+1}$}} \\ 
\bar{C}^{p,q} &
\stackrel{d}{\longrightarrow} & \bar{C}^{p,q+1}
\end{array}
\]

A homomorphism of double complexes,
\[
\phi^{\bullet,\bullet}: \Ch^{\bullet,\bullet} \longrightarrow \bar{\Ch}^{\bullet,\bullet},
\]
induces an homomorphism of the corresponding total complexes which we
will denote by,
\[
\Tot^{\bullet}(\phi^{\bullet,\bullet}): \Tot^{\bullet}(\Ch^{\bullet,\bullet}) \longrightarrow \Tot^{\bullet}(\bar{\Ch}^{\bullet,\bullet}).
\]
It also induces homomorphisms,
$'\phi_s: {'E}_s \longrightarrow {'\bar{E}}_s$
(respectively,
$''\phi_s: {''E}_s \longrightarrow {''\bar{E}}_s$)
between the associated spectral sequences (corresponding either to the
row-wise or column-wise filtrations).
For the precise definition of homomorphisms of spectral sequences,
see \cite{Mcleary}.
We will need the following useful fact 
(see \cite{Mcleary}, page 66, Theorem 3.4 for a proof).

\begin{proposition}
\label{prop:spectral}
If ${'\phi}_s$  (respectively, ${''\phi}_s$)
is an isomorphism for some $s \geq 1$, then
${'E}_r^{p,q}$ and ${'\bar{E}}_r^{p,q}$ 
(repectively, ${''E}_r^{p,q}$ and ${''\bar{E}}_r^{p,q}$ )
are isomorphic for
all $r \geq s$. In particular, the induced homomorphism,
\[
\Tot^{\bullet}(\phi^{\bullet,\bullet}): \Tot^{\bullet}(\Ch^{\bullet,\bullet}) \longrightarrow 
\Tot^{\bullet}(\bar{C}^{\bullet,\bullet})
\]
is a quasi-isomorphism.
\end{proposition}

\subsection{The Mayer-Vietoris Double Complex}
\label{subsec:MV}
Let $A_1,\ldots,A_n$ be sub-complexes of a finite simplicial complex
$A$ such that $A = A_1 \cup \cdots \cup A_n$. 
Note that the intersections of any number of the sub-complexes, $A_i$,
is again a sub-complex of $A$.
We will denote by $A_{\alpha_0\cdots \alpha_p}$ the sub-complex
$A_{\alpha_0} \cap \cdots \cap A_{\alpha_p}$.

Let $\Ch^i(A)$ denote the
$\Q$-vector space of $i$ co-chains of $A$, and 
$\Ch^{\bullet}(A),$ the complex
\[ \cdots \rightarrow \Ch^{q-1}(A) \stackrel{d}{\longrightarrow}
\Ch^q(A) \stackrel{d}{\longrightarrow} \Ch^{q+1}(A) \rightarrow \cdots
\]
where 
$d: \Ch^q(A) \rightarrow \Ch^{q+1}(A)$ 
are the usual co-boundary homomorphisms. 
More precisely, given $\omega \in  \Ch^q(A)$, and a 
$q+1$ simplex $[a_0,\ldots,a_{q+1}] \in A$,
\begin{equation}
\label{co-boundary}
d\omega([a_0,\ldots,a_{q+1}]) = \sum_{0 \leq i \leq q+1}
            (-1)^i \omega([a_0,\ldots, \hat{a_i}, \ldots,a_{q+1}])
\end{equation}
(here and everywhere else in the paper $\hat{}$ denotes
omission). Now extend $d\omega$ to a linear form on 
all of $C_{q+1}(A)$ by linearity, to
obtain an element of $\Ch^{q+1}(A).$

The connecting homomorphisms  are ``generalized'' restrictions
and are defined below.

The {\em generalized Mayer-Vietoris sequence} is the following
exact sequence of vector spaces.
\[
0 \longrightarrow \Ch^{\bullet}(A) \stackrel{r^{\bullet}}{\longrightarrow}
\bigoplus_{1 \leq \alpha_0 \leq n} \Ch^{\bullet}(A_{\alpha_0}) 
\stackrel{\delta^{0,\bullet}}{\longrightarrow}
\bigoplus_{1 \leq \alpha_0<\alpha_1 \leq n} \Ch^{\bullet}(A_{\alpha_0\cdot \alpha_1}) 
 \stackrel{\delta^{1,\bullet}}{\longrightarrow}
\cdots 
\]
\[
\bigoplus_{1 \leq \alpha_0 < \cdots < \alpha_p \leq n}\Ch^{\bullet}(A_{\alpha_0\cdots \alpha_p})
 \stackrel{\delta^{p-1,\bullet}}{\longrightarrow}
\bigoplus_{1 \leq \alpha_0< \cdots <\alpha_{p+1} \leq n}\Ch^{\bullet}(A_{\alpha_0\cdots \alpha_{p+1}})
 \stackrel{\delta^{p,\bullet}}{\longrightarrow}
\cdots
\]
where $r^{\bullet}$ is induced by restriction and the connecting homomorphisms 
$\delta^{p,\bullet}$ are as follows.

Given an $\omega \in \bigoplus_{\alpha_0< \cdots <\alpha_p}\Ch^q(A_{\alpha_0\cdots \alpha_p})$ 
we define $\delta^{p,q}(\omega)$ as follows:

First note that 
$\delta^{p,q}\omega \in \bigoplus_{\alpha_0< \cdots <\alpha_{p+1}}\Ch^q(A_{\alpha_0 \cdots \alpha_{p+1}})$, and it suffices to define
\[
(\delta^{p,q}\omega)_{\alpha_0,\ldots,\alpha_{p+1}}
\] 
for each $(p+2)$-tuple
$1 \leq \alpha_0< \cdots <\alpha_{p+1} \leq n$.
Note that, $(\delta^{p,q}\omega)_{\alpha_0,\ldots,\alpha_{p+1}}$ is a linear
form on the vector space, $C_q(A_{\alpha_0\cdots \alpha_{p+1}})$, and
hence is determined by its values on the $q$-simplices 
in the complex $A_{\alpha_0\cdots \alpha_{p+1}}$. Furthermore, 
each $q$-simplex, $s \in  A_{\alpha_0\cdots \alpha_{p+1}}$ is automatically
a simplex of the complexes 
\[
A_{\alpha_0\cdots\hat{\alpha_i}\cdots \alpha_{p+1}}, \; 0 \leq i \leq p+1.
\]

We define,
\[
\label{delta}
(\delta^{p,q} \omega)_{\alpha_0,\ldots,\alpha_{p+1}}(s) =  
 \sum_{0 \leq j \leq p+1} (-1)^i \omega_{\alpha_0,\ldots,\hat{\alpha_j},\ldots,\alpha_{p+1}}
(s).
\]
The fact that the generalized Mayer-Vietoris  sequence
is exact is classical (see \cite{Rotman} or \cite{B03} for example).

We now define the Mayer-Vietoris double complex of the complex $A$
with respect to the subcomplexes $A_{\alpha_0}, 1 \leq \alpha_0 \leq n$, 
which we will denote by  $\N^{\bullet,\bullet}(A)$ (we suppress the
dependence of the complex on sub-complexes $A_{\alpha_0}$ in the notation
since this dependence will be clear from context). 

\begin{definition}
\label{def:MV}
The Mayer-Vietoris double complex of a simplicial complex $A$
with respect to the subcomplexes $A_{\alpha_0}, 1 \leq \alpha_0 \leq n$, 
$\N^{\bullet,\bullet}(A)$, is the double complex defined by,
$$
\displaylines{
\N^{p,q}(A) = 
\bigoplus_{1 \leq \alpha_0< \cdots <\alpha_p \leq n}
\Ch^q(A_{\alpha_0\cdots\alpha_p}).
}
$$
The horizontal differentials are as defined above. The vertical
differentials are those induced by the ones in the different
complexes, $\Ch^{\bullet}(A_{\alpha_0\cdots\alpha_p}).$

$\N^{\bullet,\bullet}(A)$ is depicted in the following figure.
\begin{diagram}
\uTo && \uTo &&& \\
\bigoplus_{\alpha_0} \Ch^2(A_{\alpha_0}) & \rTo &
\bigoplus_{\alpha_0 < \alpha_1} \Ch^2(A_{\alpha_0\cdot\alpha_1}) & \rTo & \ldots \\
\uTo && \uTo &&& \\
\bigoplus_{\alpha_0} \Ch^1(A_{\alpha_0}) & \rTo &
\bigoplus_{\alpha_0 < \alpha_1} \Ch^1(A_{\alpha_0\cdot\alpha_1}) & \rTo & \ldots \\
\uTo && \uTo &&& \\
\bigoplus_{\alpha_0} \Ch^0(A_{\alpha_0}) & \rTo &
\bigoplus_{\alpha_0 < \alpha_1} \Ch^0(A_{\alpha_0\cdot\alpha_1}) & \rTo & \ldots \\
\end{diagram}

For any $t \geq 0,$
we denote by $\N_t^{\bullet,\bullet}(A)$ the following truncated complex.
\[
\begin{array}{cccc}
\N_t^{p,q}(A) & = & \N^{p,q}(A), & 
0 \leq p+q \leq t, \cr
\N_t^{p,q}(A) & = &  0, & \mbox{otherwise}.\cr
\end{array}
\]
\end{definition}

The following proposition is classical (see \cite {Rotman} or \cite{B03} 
for a proof) and 
follows from the exactness of the generalized Mayer-Vietoris sequence.

\begin{proposition}
\label{prop:MV}
The  spectral sequences,
${'E}_r,{''E}_r$,
associated to $\N^{\bullet,\bullet}(A)$ converge to $\HH^*(A)$ and thus,
$$
\displaylines{
\HH^*(\Tot^{\bullet}(\N^{\bullet,\bullet}(A))) \cong \HH^*(A).
}
$$
Moreover, the homomorphism 
\[
\psi^{\bullet}: \Ch^{\bullet}(A) \rightarrow \Tot^{\bullet}(\N^{\bullet,\bullet}(A))
\]
induced by the homomorphism $r^{\bullet}$ (in the generalized Mayer-Vietoris
sequence) is a quasi-isomorphism.
\end{proposition}

We denote by  $\Ch^{\bullet}_{\ell+1}(A)$ the truncation of the
complex $\Ch^{\bullet}(A)$ after the $(\ell+1)$-st term.
As an immediate corollary we have that,
\begin{corollary}
\label{cor:MV}
For any $\ell \geq 0$,
the homomorphism 
\begin{equation}
\label{eqn:restriction}
\psi_{\ell+1}^{\bullet}: \Ch^{\bullet}_{\ell+1}(A) \rightarrow \Tot^{\bullet}(\N_{\ell+1}^{\bullet,\bullet}(A))
\end{equation}

induced by the homomorphism $r^{\bullet}$ (in the generalized Mayer-Vietoris
sequence) is a quasi-isomorphism.
Hence, for $0 \leq i \leq \ell,$
$$
\displaylines{
\HH^i(\Tot^{\bullet}(\N_{\ell+1}^{\bullet,\bullet}(A))) \cong \HH^i(A).
}
$$
\end{corollary}

\begin{remark}
\label{rem:induction}
Notice that in the truncated Mayer-Vietoris double complex,
$\N_t^{\bullet,\bullet}(A)$, the $0$-th column is a complex
having at most $t+1$ non-zero terms, the first column can have at most
$t$ non-zero terms, and in general the $i$-th column has at most $t+1 - i$
non-zero terms. This observation plays a crucial role in the inductive
argument used later in the paper (in the proof of Proposition 
\ref{prop:bound}).
\end{remark}

\section{Double complexes associated to certain covers}
\label{sec:covering}
We begin with a definition.
\begin{definition}
\label{def:P-closed}
Let ${\mathcal P}$ be a finite subset of $\R[X_1,\ldots,X_k]$.
A ${\mathcal P}$-closed formula
is a formula  constructed as follows:
\begin{description}
\item
For each $P \in {\mathcal P}$,
$$P=0, P \ge 0, P\le 0,$$
\noindent are ${\mathcal P}$-closed formulas.
\item
If $\Phi_1$ and $\Phi_2$ are
${\mathcal P}$-closed formulas, $\Phi_1 \wedge \Phi_2$ and
$\Phi_1 \vee \Phi_2$
 are ${\mathcal P}$-closed formulas.
\end{description}

Clearly, $\RR(\Phi) = \{ x \subset \R^k \mid  \Phi(x)\}$, the
realization of a ${\mathcal P}$-closed formula $\Phi,$ is a closed
semi-algebraic set and we call such a set a ${\mathcal P}$-closed 
semi-algebraic set.
\end{definition}

In this section, we consider a fixed family of polynomials,
${\mathcal P} \subset \R[X_1,\ldots,X_k],$ as well as a fixed  
${\mathcal P}$-closed and bounded semi-algebraic set, $S \subset \R^k$. 
We also fix a number, $\ell, 0 \leq \ell \leq k.$

We define below (in Section \ref{subsec:defineadmissible})
a finite set of indices, $\A_S$, which we call the set of admissible
indices, and a map that associates to each $\alpha \in \A_S$ a closed and
bounded semi-algebraic subset $X_\alpha \subset S$, which we call an
admissible subset.
To each $\alpha  \in \A_S$, 
we associate its level, denoted $\level(\alpha)$, which is an integer
between $0$ and $\ell$. The set $\A_S$ will be partially ordered, and we
denote by $\ancestor(\alpha)\subset \A_S$, the set of ancestors
of $\alpha$ under this partial order. For $\alpha,\beta \in \A_S$,
$\beta \in \ancestor(\alpha),$ implies that  
$X_\alpha \subset X_\beta$.

For each admissible index $\alpha \in \A_S$, 
we define  a double complex, $\M^{\bullet,\bullet}(\alpha)$, 
such that 
\[
\HH^i(\Tot^{\bullet}(\M^{\bullet,\bullet}(\alpha))) \cong \HH^i(X_\alpha),\;
0 \leq i \leq \ell - \level(\alpha).
\]
 
The main idea behind the construction of the double complex 
$\M^{\bullet,\bullet}(\alpha)$ is as follows. 
Associated to  any cover of
$X_\alpha$ there exists  a double complex (the Mayer-Vietoris double complex)
arising from the generalized Mayer-Vietoris exact sequence (see \cite{B03}).
If the individual sets of the cover of $X$ are all contractible, then the
first column of the Mayer-Vietoris double complex is  zero except at
the first row.
The cohomology groups of the associated total complex of the Mayer-Vietoris
double complex are isomorphic to those of $X_\alpha$ 
and thus in order to compute
$b_0(X_\alpha),\ldots,b_{\ell - \level(\alpha)}(X_\alpha),$
it suffices to compute a suitable
truncation of the Mayer-Vietoris double complex. However, computing 
(even the truncated) Mayer-Vietoris double complex directly within a 
singly exponential time complexity
is not possible by any known method, since we are unable 
to compute triangulations of semi-algebraic sets in singly exponential time.
However, 
making use of the cover construction recursively, we are
able to compute another double complex, 
$\M^{\bullet,\bullet}(\alpha)$, 
which has much smaller size but whose associated
total complex  is  quasi-isomorphic to the truncated
Mayer-Vietoris double complex and hence has isomorphic cohomology groups
(see Proposition \ref{prop:main} below).
The construction of 
$\M^{\bullet,\bullet}(\alpha)$ is possible in singly exponential time since the
covers can be computed in singly exponential time.

Finally, given any closed and bounded semi-algebraic set $X \subset \R^k$,
we will denote by ${\mathcal C}'(X)$, 
a fixed cover of $X$
(we will use the construction in \cite{BPR04} to compute
such a cover).

\subsection{Admissible sets and Covers}
\label{subsec:defineadmissible}
We now define $\A_S$,
and for each $\alpha \in \A_S$ a cover
${\mathcal C}(\alpha)$ of $X_\alpha$ 
obtained by enlarging the cover ${\mathcal C}'(X_\alpha)$. 

\begin{definition}(Admissible indices and covers)
\label{def:admissible}
We define $\A_S$ by induction on level.
\begin{enumerate}
\item
Firstly, 
$0 \in \A_S$, $\level(0) = 0,$ $X_0 = S$,
and ${\mathcal C}(0) = {\mathcal C}'(S).$ 
The admissible indices at level $1$ consists of all formal
products,
$\beta = \alpha_0\cdot\alpha_1\cdots\alpha_{j-1}\cdot\alpha_j$,
with $\alpha_i \in {\mathcal C}(0)$ and $0 \leq j \leq \ell+1$, 
and we define the associated semi-algebraic set by, 
\[
X_\beta = X_{\alpha_0} \cap \cdots \cap X_{\alpha_j}.
\]  
For each 
$\{\alpha_0,\ldots,\alpha_m\} \subset \{\beta_0,\ldots,\beta_n \} 
\subset {\mathcal C}(0)$,  
with $n \leq \ell+1$, 
\[
\alpha_0\cdots\alpha_m \in  \ancestor(\beta_0\cdots\beta_n),
\] 
and 
$0 \in \ancestor(\beta_0\cdots\beta_n)$.

\item 
We now inductively define the admissible indices  at level ${i+1}$, in terms
of the admissible indices at level $\leq i$.
For each $\alpha \in \A_S$ at level $i$, we define
${\mathcal C}(\alpha)$ as follows. 
Let $\ancestor(\alpha) = \{\alpha_1,\ldots,\alpha_N\}$. 
Then,
$$
\displaylines{
{\mathcal C}(\alpha) = \dot{\bigcup}_{\beta_i \in {\mathcal C}(\alpha_i), 1 \leq i \leq N}
{\mathcal C}'(\beta_1  \cdots \beta_N \cdot \alpha),
}
$$
where $\dot\bigcup$ denotes the disjoint union.
All formal products, $\beta = \alpha_0\cdot\alpha_1\cdots\alpha_j$,
with $\alpha_i \in {\mathcal C}(\alpha)$ and $0 \leq j \leq \ell-i+1$ 
are in $\A_S$, and we define 
\[
X_\beta = X_{\alpha_0} \cap \cdots \cap X_{\alpha_j},
\]  
and $\level(\beta) = i+1.$

For each 
$\{\alpha_0,\ldots,\alpha_m\} \subset \{\beta_0,
\ldots,\beta_n \} \subset {\mathcal C}(\alpha)$,  
with $n \leq \ell-i+1$, 
\[
\alpha_0\cdots\alpha_m \in  \ancestor(\beta_0\cdots\beta_n),
\] 

and $\alpha \in \ancestor(\beta_0\cdots\beta_n).$

Moreover, for $\alpha' \in {\mathcal C}'(\beta_1 \cdot \cdots\cdot \beta_N 
\cdot \alpha)$, each
$\beta_i$ is an ancestor of $\alpha'$.
We transitively close the ancestor relation, so that ancestor of 
an ancestor is also an ancestor.
Moreover, if 
$\alpha_0\cdots\alpha_m, \beta_0\cdots\beta_n \in \A_S$ are such that
for every  $j \in \{1,\ldots, n\}$ there exists $i \in \{1,\ldots, m\}$ 
such that
$\alpha_i$ is an ancestor of $\beta_j$, then $\alpha_0\cdots\alpha_m$ is
an ancestor of $\beta_0\cdots\beta_n$.

Finally, the set of admissible indices at level $i+1$ is
\[
\dot\bigcup_{\alpha \in \A_S, \level(\alpha) = i}
\{\alpha_0\cdot\alpha_1\cdots\alpha_j \mid 
\alpha_i \in {\mathcal C}(\alpha), 0 \leq j \leq \ell-i+1\}
\]
\end{enumerate}
\end{definition}

Observe that by the above definition, if $\alpha,\beta \in \A_S$
and $\beta \in \ancestor(\alpha)$,
then each $\alpha' \in {\mathcal C}(\alpha)$ has a unique ancestor 
in  each  ${\mathcal C}(\beta)$, which we will denote by
$a_{\alpha,\beta}(\alpha')$ and the mapping,
$a_{\alpha,\beta}: {\mathcal C}(\alpha) \rightarrow {\mathcal C}(\beta)$
is injective. 

\hide{
Now, suppose that we have a procedure for computing  ${\mathcal C}'(X)$, 
for any given ${\mathcal P}'$-closed and bounded semi-algebraic set, $X$,
where $\#{\mathcal P}' = m$ and $\deg(P) \leq D,$ for $P \in {\mathcal P}' $.
Moreover, suppose that
the number and the degrees of the polynomials appearing in the
output of this procedure is bounded by 
$(mD)^{k^{O(1)}}.$
}
Now, suppose that we have a procedure for computing  ${\mathcal C}'(X)$, 
for any given ${\mathcal P}'$-closed and bounded semi-algebraic set, $X$,
such that
the number and the degrees of the polynomials appearing the descriptions
of the semi-algebraic sets, $X_\alpha, \alpha \in {\mathcal C}'(X)$,
is bounded by 
\begin{equation}
\label{eqn:c_1}
D^{k^{c_1}},
\end{equation}
where $c_1 > 0$ is some absolute constant, and 
$D = \sum_{P \in {\mathcal P}'} \deg(P).$
  
Using  the above procedure for computing ${\mathcal C}'(X)$, and
the definition of $\A_S$,  we have the following
quantitative bounds on $\#\A_S$  and the 
 semi-algebraic sets $X_\alpha, \alpha \in \A_S$,
which is crucial in proving the complexity bound of our algorithm.

\begin{proposition}
\label{prop:bound}
Let $S \subset \R^k$ be a ${\mathcal P}$-closed semi-algebraic set, where
${\mathcal P} \subset \R[X_1,\ldots,X_k]$ is a family of $s$ polynomials of 
degree at most $d$. Then $\#A_S$, as well as 
the number of polynomials used to define the semi-algebraic sets 
$X_\alpha, \alpha \in \A_S$ and the  
the degrees of these polynomials, are all bounded by
$(sd)^{k^{O(\ell)}}.$
\end{proposition}

\begin{proof}{Proof:}
Given $\alpha \in \A_S$ with $\level(\alpha) = j$, 
we first prove by induction on $\level(\alpha)$ that
\[
\#\ancestor(\alpha) \leq 2^{(j+1)(\ell+3) - j(j+1)/2}.
\]
The claim is clearly true if $\level(\alpha) = 0$.
Otherwise, from the definition of $\A_S$, 
there exists $\beta \in \A_S$,
with $\level(\beta) = j-1$, such that 
$\alpha  = \gamma_0\cdots\gamma_m,$
$\gamma_i \in {\mathcal C}(\beta)$ and $m \leq \ell -j +2 .$

For each $\gamma_i$, we have
\[
\ancestor(\gamma_i) = \ancestor(\beta) \cup \{ a_{\beta,\theta}(\gamma_i) 
\;\mid\; 
\theta  \in \ancestor(\beta)\},
\]

and it follows that,
$$
\displaylines{
\ancestor(\alpha) = \ancestor(\beta) \cup 
\{a_{\beta,\theta}(\gamma_{i_0})\cdots a_{\beta,\theta}(\gamma_{i_n})
\;\mid\; \cr
\theta  \in \ancestor(\beta),
\{i_0,\ldots,i_n\} \subset \{1,\ldots,m\}
\}.
}
$$

Hence, 
$$
\begin{array}{ccc}
\#\ancestor(\alpha) &  =  &  \#\ancestor(\beta)\cdot 2^{m}\\
                    &\leq &  \#\ancestor(\beta)\cdot 2^{\ell - j + 3}\\
                    &\leq &  2^{\sum_{i=0}^{j} (\ell - i + 3)} \\
                    &  =  &  2^{(j+1)(\ell+3) - j(j+1)/2} \\
                    &\leq &  2^{c_2 \ell^2},
\end{array}
$$
for some absolute constant $c_2$.

We now prove again by induction on the level 
that there exists an absolute constant
$c > 0$, such that the  number of elements of $\A_S$ of level $ \leq j$,
as well as the number of polynomials needed to define the associated
semi-algebraic sets, and 
the degrees of these polynomials, are all bounded by
$(sd)^{k^{c j}}.$ 

The claim is clear for level $0$. 
Now assume that the claim holds for level $< j$.
As before, given $\alpha \in \A_S$ with
$\level(\alpha) = j$, there exists $\beta \in \A_S$
with $\level(\beta) = j-1$, such that
$\alpha  = \gamma_0\cdots\gamma_m,$
$\gamma_i \in {\mathcal C}(\beta)$ and $m \leq \ell -j +2 .$
We have that $\#\ancestor(\beta) \leq 2^{c_2 \ell^2}$ by the previous
paragraph. 
Let $\ancestor(\beta) = \{ \theta_1,\ldots,\theta_N\}$. Then,
\[
\#{\mathcal C}(\theta_i) \leq (sd)^{k^{c (j-1)}},
\]
for $1 \leq i \leq N$ by the induction hypothesis.

In order to bound,
$$
\displaylines{
\#{\mathcal C}(\beta) = \#\bigcup_{\beta_i \in {\mathcal C}(\theta_i), 
1 \leq i \leq N}
{\mathcal C}'(\beta_1 \cdot \cdots\cdot  \beta_N \cdot \beta),
}
$$
first observe that 
$N \leq 2^{c_2 \ell^2}$ and hence
the union on the right hand side is over
an index set of cardinality bounded by, 
\[
(sd)^{k^{c (j-1)} 2^{c_2\ell^2}},
\]
and each set in the union has cardinality bounded by,

$$
\begin{array}{ccc}
M  &=& (2^{c_2\ell^2}(sd)^{k^{c (j-1)}})^{k^{c_1}} \\
   &=& 2^{c_2\ell^2 k^{c_1}} (sd)^{k^{cj - (c - c_1)}} ,
\end{array}
$$
where $c_1$ is the constant defined before  in (\ref{eqn:c_1}) above.

Thus, the total number of admissible indices at level $j$ is bounded
by the total number of admissible indices at level $j-1$ times
$\sum_{0 \leq i \leq \ell -j+3} {M \choose i}$.

It follows that if $c$ chosen large enough with respect to the 
constants $c_1,c_2$, then for all $k$ large enough, the 
the total number of admissible indices at level $j$ is at most,
\[
(sd)^{k^{c j}}. 
\]

The bounds on the number and degrees of polynomials appearing in the
description can be proved similarly using the same induction scheme.
\end{proof}

\subsection{Double Complex Associated to a Cover}
Given the different covers described above, 
we now associate to each 
$ \alpha \in \A_S$  
a double complex,
${\mathcal M}^{\bullet,\bullet}(\alpha),$ 
and for  every $\beta \in \A_S$, such that $\alpha \in \ancestor(\beta)$,
and $\level(\alpha) = \level(\beta)$,
a restriction homomorphism:
\[
r_{\alpha,\beta}^{\bullet,\bullet}:  \M^{\bullet,\bullet}(\alpha) 
\rightarrow \M^{\bullet,\bullet}(\beta),
\]
satisfying the following:

\begin{enumerate}
\item
\begin{equation}
\label{eqn:iso}
\HH^i(\Tot^{\bullet}(\M^{\bullet,\bullet}(\alpha)))  \cong 
\HH^i(X_\alpha), \;\mbox{\rm for} \;
0 \leq i \leq \ell - \level(\alpha).
\end{equation}
\item
The restriction homomorphism
\[
r_{\alpha,\beta}^{\bullet,\bullet}:  \M^{\bullet,\bullet}(\alpha) 
\rightarrow \M^{\bullet,\bullet}(\beta),
\]
induces the restriction homomorphisms between the cohomology
groups:
$$
\displaylines{
r^*_{\alpha,\beta}: \HH^i(X_\alpha) \rightarrow \HH^i(X_\beta)
}
$$ 
for $0 \leq i \leq \ell - \level(\alpha)$  via the isomorphisms in (\ref{eqn:iso}).
\end{enumerate}

We now describe the construction of the double complex $\M^{\bullet,\bullet}(\alpha)$ and
prove that it has the properties stated above. The 
double complex $\M^{\bullet,\bullet}(\alpha)$, is constructed inductively using induction on $\level(\alpha)$.

\begin{definition}
\label{def:double}
The base case is when 
$\level(\alpha) = \ell$.
In this case the double complex, $\M^{\bullet,\bullet}(\alpha)$
is defined by:
\[
\begin{array}{ccll}
\M^{0,0}(\alpha) & = & \bigoplus_{{\alpha_0}\; \in\; {\mathcal C}(\alpha)}\; \HH^0(X_{\alpha_0}),
\cr
\M^{1,0}(\alpha) & = & \bigoplus_{{\alpha_0}, {\alpha_1} \;\in\; {\mathcal C}(\alpha)}\; \HH^0(X_{\alpha_0\cdot\alpha_1}), \cr
\M^{p,q}(\alpha) & = &  0, \; \mbox{if} \; q > 0 \;\mbox{or} \; p > 1 . \cr
\end{array}
\]

This is shown diagramatically below.
\begin{diagram}
0&\rTo&0&\rTo&0& \\
\uTo&&\uTo&&\uTo&\\
0 & \rTo& 0 &\rTo&0&\\
\uTo &&\uTo &&\uTo&\\
\bigoplus_{{\alpha_0} \in {\mathcal C}(\alpha)} 
\HH^0(X_{\alpha_0}) & \rTo^\delta & 
\bigoplus_{{\alpha_0}, {\alpha_1} \in {\mathcal C}(\alpha)} 
\HH^0(X_{\alpha_0\cdot \alpha_1}) & \rTo & 0&
\end{diagram}

The only non-trivial homomorphism in the above complex, 
$$
\displaylines{
\delta:\bigoplus_{{\alpha_0} \in {\mathcal C}(\alpha)} 
\HH^0(X_{\alpha_0}) 
\longrightarrow
\bigoplus_{{\alpha_0},{\alpha_1} \in {\mathcal C}(\alpha)} 
\HH^0(X_{\alpha_0\cdot \alpha_1}) 
}
$$
is defined as follows.

$\delta(\phi)_{\alpha_0,\alpha_1} = 
(\phi_{\alpha_1} - \phi_{\alpha_0})|_{X_{\alpha_0\cdot\alpha_1}}$ for
$\phi \in \bigoplus_{{\alpha_0} \in {\mathcal C}(\alpha)} \HH^0(X_{\alpha_0}).$

For every $\beta \in \A_S$, such that $\alpha \in \ancestor(\beta)$,
and $\level(\alpha) = \level(\beta) = \ell$,
we define $r_{\alpha,\beta}^{0,0}: \M^{0,0}(\alpha) \rightarrow \M^{0,0}(\beta),$
as follows.

Recall that,
$\displaystyle{
\M^{0,0}(\alpha) = \bigoplus_{\alpha_0 \;\in\; {\mathcal C}(\alpha)} 
\; \HH^0(X_{\alpha_0}),
}
$
and 
$
\displaystyle{
\M^{0,0}(\beta) = \bigoplus_{\beta_0 \;\in\; {\mathcal C}(\beta)} \; \HH^0(X_{\beta_0}).
}
$

For $\phi \in \M^{0,0}(\alpha)$ and $\beta_0 \in {\mathcal C}(\beta)$ 
we define,
$$
\displaylines{
r_{\alpha,\beta}^{0,0}(\phi)_{\beta_0} = \phi_{a_{\beta,\alpha}(\beta_0)}
|_{X_{\beta_0}}.
}
$$

We define
$r_{\alpha,\beta}^{1,0}: \M^{1,0}(\alpha) \rightarrow \M^{1,0}(\beta),$
in a similar manner. More precisely,
for $\phi \in \M^{0,0}(\alpha)$ and
$\beta_0,\beta_1 \in {\mathcal C}(\beta)$,
we define
$$
\displaylines{
r_{\alpha,\beta}^{1,0}(\phi)_{\beta_0,\beta_1} = 
\phi_{a_{\beta,\alpha}(\beta_0)\cdot a_{\beta,\alpha}(\beta_1)}|_{X_{\beta_0\cdot\beta_1}}.
}
$$ 


(The inductive step)
In general the $\M^{p,q}(\alpha)$ are defined as follows using induction
on $\level(\alpha)$ and with $n_\alpha= \ell - \level(\alpha)+1.$
$$
\displaylines{
\begin{array}{ll}
\M^{0,0}(\alpha) =  \bigoplus_{{\alpha_0} \;\in\; {\mathcal C}(\alpha)} 
\;\HH^0(X_{\alpha_0}),& \cr
\M^{0,q}(X) =   0, & 0< q, \cr
\M^{p,q}(\alpha) =   \bigoplus_{\alpha_0< \cdots <\alpha_p, \;
{\alpha_i}  \in {\mathcal C}(\alpha)}
\; \Tot^q(\M^{\bullet,\bullet}({\alpha_0 \cdots \alpha_p})),  & 0 < p,\;
0< p+q \leq n_\alpha,  \cr
\M^{p,q}(\alpha) =   0, &\;\mbox{else}.
\end{array}
}
$$

The double complex $\M^{\bullet,\bullet}(\alpha)$ is shown  in the following 
diagram:
{\tiny
\begin{diagram}
0 & \rTo & 0 & \rTo & 0 & \cdots & 0\\
\uTo && \uTo && \uTo && \uTo  \\
0 & \rTo &
\bigoplus_{\alpha_0 < \alpha_1}\Tot^{n_\alpha-1}(\M^{\bullet,\bullet}({\alpha_0\cdot \alpha_1})) 
&\rTo^{\delta}&
0&\cdots&0 \\
\uTo_{d} &&\uTo_{d}&&\uTo_{d} &&\uTo_{d} \\
0 & \rTo &
\bigoplus_{\alpha_0 < \alpha_1}\Tot^{n_\alpha-2}(\M^{\bullet,\bullet}({\alpha_0\cdot\alpha_1})) 
&\rTo^{\delta}&
\bigoplus_{\alpha_0 < \alpha_1 < \alpha_2}\Tot^{n_\alpha-2}(\M^{\bullet,\bullet}
({\alpha_0\cdot\alpha_1\cdot \alpha_2}))
&\cdots &0 \\
\vdots &\vdots& \vdots&\vdots& \vdots&\vdots&\vdots\\
0 & \rTo &
\bigoplus_{\alpha_0 < \alpha_1}\Tot^2(\M^{\bullet,\bullet}({\alpha_0\cdot\alpha_1})) &
\rTo^{\delta}&
\bigoplus_{\alpha_0< \alpha_1 < \alpha_2}
\Tot^2(\M^{\bullet,\bullet}({\alpha_0\cdot\alpha_1\cdot\alpha_2}))
&\cdots&0 \\
\uTo_{d} &&\uTo_{d}&&\uTo_{d} && \uTo_{d}\\
0 & \rTo &
\bigoplus_{\alpha_0 < \alpha_1}\Tot^1(\M^{\bullet,\bullet}({\alpha_0\cdot\alpha_1})) 
&\rTo^{\delta}&
\bigoplus_{\alpha_0 < \alpha_1 < \alpha_2}\Tot^1(\M^{\bullet,\bullet}({\alpha_0\cdot\alpha_1\cdot\alpha_2}))
&\cdots&0 \\
\uTo_{d} &&\uTo_{d}&&\uTo_{d} &&\uTo_{d} \\
\bigoplus_{X_{\alpha_0} \in {\mathcal C}_X} \HH^0(X_{\alpha_0}) & \rTo^{\delta} & 
\bigoplus_{\alpha_0 < \alpha_1}\Tot^0(\M^{\bullet,\bullet}({\alpha_0\cdot\alpha_1})) &\rTo^{\delta}&
\bigoplus_{\alpha_0 < \alpha_1 < \alpha_2}\Tot^0(\M^{\bullet,\bullet}({\alpha_0\cdot\alpha_1\cdot \alpha_2})) &\cdots&
\bigoplus_{\alpha_0< \cdots <\alpha_{n_\alpha}}\Tot^0(\M^{\bullet,\bullet}
({\alpha_0\cdots\alpha_{n_\alpha}}))
\end{diagram}
}

The vertical homomorphisms, $d$, in $\M^{\bullet,\bullet}(\alpha)$ are those 
induced by the differentials in the various 
$$
\Tot^{\bullet}(\M^{\bullet,\bullet}({\alpha_0\cdots\alpha_p})),
{\alpha_i} \in {\mathcal C}(\alpha).
$$

The horizontal ones are defined by generalized restriction as follows.
Let 
\[\phi \in \bigoplus_{\alpha_0< \cdots <\alpha_p, \alpha_i \in 
{\mathcal C}(\alpha)} 
\Tot^q(\M^{\bullet,\bullet}({\alpha_0\cdots\alpha_p})),
\]

with 
\[
\phi_{\alpha_0,\ldots,\alpha_p} = \bigoplus_{0 \leq j \leq q}
\phi_{\alpha_0,\ldots,\alpha_p}^j,
\]
and 
\[
\phi_{\alpha_0,\ldots,\alpha_p}^j \in \M^{j,q-j}({\alpha_0\cdots\alpha_p}).
\]

We define,
$$\displaylines{
\delta:
\bigoplus_{\alpha_0< \cdots <\alpha_p, \alpha_i \in {\mathcal C}(\alpha)} 
\Tot^q(\M^{\bullet,\bullet}({\alpha_0\cdots\alpha_p})) \longrightarrow
\bigoplus_{\alpha_0< \cdots <\alpha_{p+1}}
\Tot^q(\M^{\bullet,\bullet}({\alpha_0\cdots\alpha_{p+1}}))
}
$$
by
$$
\displaylines{
\delta(\phi)_{\alpha_0,\ldots, \alpha_{p+1}} = 
\bigoplus_{0 \leq i \leq p+1} (-1)^i 
\bigoplus_{0 \leq j \leq q}
r_{{\alpha_0\cdots\hat{\alpha_i}
\cdots\alpha_{p+1}},{\alpha_0\cdots\alpha_{p+1}}}^{j,q-j}
(\phi_{\alpha_0,\ldots,\hat{\alpha_i},
\ldots,\alpha_{p+1}}^j),
}
$$
noting that for each $i, 0 \leq i \leq p+1,$ 
${\alpha_0\cdots\hat{\alpha_i}\cdots\alpha_{p+1}}$
is an ancestor of ${\alpha_0\cdots\alpha_{p+1}}$, and

\[
\level({\alpha_0\cdots\hat{\alpha_i}\cdots\alpha_{p+1}})=
\level({\alpha_0\cdots\alpha_{p+1}}) = \level(\alpha)+ 1,
\]
and hence the homomorphisms 
$r_{{\alpha_0\cdots\hat{\alpha_i}
\cdots\alpha_{p+1}},{\alpha_0\cdots\alpha_{p+1}}}^{j,q-j}$
are already defined by induction.

Now let,
$\alpha,\beta \in \A_S$ 
with $\alpha$ an ancestor of $\beta$ and
$\level(\alpha) = \level(\beta).$
We define the restriction homomorphism, 
$$
r_{\alpha,\beta}^{\bullet,\bullet}:  
\M^{\bullet,\bullet}(\alpha) \longrightarrow \M^{\bullet,\bullet}(\beta)
$$
as follows.

As before,
for $\phi \in \M^{0,0}(\alpha)$ and 
$\beta_0 \in {\mathcal C}(\beta)$ we define,
$$
\displaylines{
r_{\alpha,\beta}^{0,0}(\phi)_{\beta_0} = \phi_{a_{\beta,\alpha}(\beta_0)}|_{X_{\beta_0}}.
}
$$
For $0 < p, 0< p+q \leq \ell-\level(\alpha)+1,$
we define
\[
r_{\alpha,\beta}^{p,q}: \M^{p,q}(\alpha) \rightarrow \M^{p,q}(\beta),
\]
as follows.

Let 
$
\displaystyle{
\phi \in \M^{p,q}(\alpha) = \bigoplus_{\alpha_0< \cdots <\alpha_p, \;
{\alpha_i}  \in {\mathcal C}(\alpha)}
\; \Tot^q(\M^{\bullet,\bullet}({\alpha_0\cdots\alpha_p})).
}
$ 
We define,
\[
r_{\alpha,\beta}^{p,q}(\phi) =
\bigoplus_{\beta_0< \cdots <\beta_p, \beta_i \in {\mathcal C}(\beta)} 
\bigoplus_{0 \leq i \leq  q}r_{a_{\beta,\alpha}({\beta_0\cdots\beta_p}),
{\beta_0\cdots\beta_p}}^{i,q-i} \phi_{a_{\beta,\alpha}
({\beta_0}),\ldots, a_{\beta,\alpha}({\beta_p})}^i,
\]
where 
$a_{\beta,\alpha}({\beta_0\cdots\beta_p}) = 
a_{\beta,\alpha}({\beta_0})\cdots a_{\beta,\alpha}({\beta_p}).$
Note that, each $a_{\beta,\alpha}({\beta_i}), 0 \leq i \leq p$ are all
distinct and belong to ${\mathcal C}(\alpha)$.
Moreover,
$$
\level(a_{\beta,\alpha}({\beta_0\cdots\beta_p}))= 
\level({\beta_0\cdots\beta_p})
= \level(\alpha) +1,
$$
and hence we can assume that the homomorphisms
$r_{a_{\beta,\alpha}({\beta_0\cdots\beta_p}),{\beta_0\cdots\beta_p}}^
{\bullet,\bullet}$ 
used in the definition of
$r_{\alpha,\beta}^{\bullet,\bullet}$ are already defined by induction.
\end{definition}

It is  easy to verify by induction on 
$\level(\alpha)$ that, $\M^{\bullet,\bullet}(\alpha)$ defined
as above, is indeed a double complex, that is the homomorphisms 
$d$ and $\delta$ satisfy  the equations,
$$ 
\displaylines{
d^2= \delta^2 = 0, 
d\; \delta + \delta\; d = 0.
}
$$


\subsection{Example}
Before proving the main properties of the complexes 
$\M^{\bullet,\bullet}(\alpha)$ defined above, 
we illustrate their construction by means of a simple example. 
We take for the set $S$, the unit sphere
$S^2 \subset \R^3$.
Even though this example looks very simple, it is actually illustrative
of the main topological ideas behind the construction of the complex
$\M^{\bullet,\bullet}(S)$ starting from a cover of $S$ by two
closed hemispheres meeting at the equator. Since the intersection of
the two hemisphere is a topological circle which is not 
contractible, Theorem \ref{the:nerve} is not applicable.
Using Theorem  \ref{the:bettione} we can compute $\HH^0(S),\HH^1(S)$, but it
is not enough to compute $\HH^2(S)$.  
The recursive
construction of ${\mathcal M}^{\bullet,\bullet}$ described in the last section
overcomes this problem and this is illustrated in the example.

\begin{figure}[hbt] 
\begin{center}
\includegraphics[width=5cm]{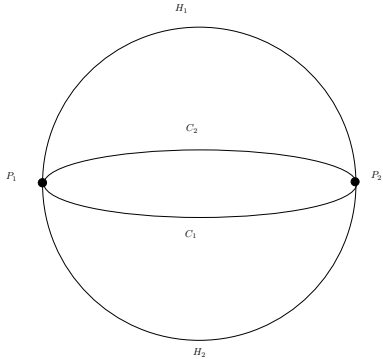}%
\end{center}
\caption{Example of $S^2 \subset \R^3$}
\label{fig:example}
\end{figure}

\begin{example}
We first fix some notations (see Figure 2).
Let $H_1$ and $H_2$ denote the
closed upper and lower hemispheres respectively. Let $H_{12} = H_1 \cap H_2$
denote the equator, 
and let $H_{12} = C_1 \cup C_2$, where $C_1,C_2$ are closed
semi-circular arcs. Finally, let $C_{12} = C_1 \cap C_2 = \{P_1,P_2\}$, where 
$P_1,P_2$ are two antipodal points.

For the purpose of this example, we will take for the 
covers ${\mathcal C'}$ the obvious ones, namely:
$$
\begin{array}{cccc}
{\mathcal C}'(S) &=& \{H_1,H_2\},& \\
{\mathcal C}'(H_i) &=& \{H_i\}, &i =1,2, \\
{\mathcal C}'(H_{12}) &=& \{C_1,C_2\},& \\
{\mathcal C}'(C_i) &=& \{C_i\},& i =1,2, \\
{\mathcal C}'(C_{12}) &=& \{P_1,P_2\}, \\
{\mathcal C}'(P_i) &=&\{P_i\},& i =1,2. \\
\end{array}
$$
Note that, in order not to complicate notations further, 
we are using the same names
for the elements of ${\mathcal C}'(\cdot)$, as well as their
associated sets. 
Strictly speaking, we should have defined,
\[
{\mathcal C}'(S) =  \{\alpha_1,\alpha_2\}, \; X_{\alpha_1} = H_1, X_{\alpha_2}
 = H_2, \ldots .\]
However, since each set occurs at most once, this
does not create confusion in this example.

Note that the elements of the sets occurring on the right are all closed, 
bounded contractible subsets of $S$.
It is now easy to check from Definition \ref{def:admissible},
that the elements of $\A_S$ in order of their levels
as follows.
\begin{enumerate}
\item Level $0$:
\[
0 \in A_S, \level(0) = 0,\]
and
$$
\begin{array}{cccc}
{\mathcal C}(0) &=& \{\alpha_1,\alpha_2\}, & X_{\alpha_1} = H_1, X_{\alpha_2} = H_2.\\
\end{array}
$$
 
\item Level $1$:
The elements of level $1$ are 
\[
\alpha_1,\alpha_2,\alpha_1\cdot\alpha_2,
\]
and
$$
\begin{array}{cccc}
{\mathcal C}(\alpha_1) &=& \{\beta_1\},  & X_{\beta_1} = H_1,\\  
{\mathcal C}(\alpha_2) &=& \{\beta_2\},  & X_{\beta_2} = H_2,\\  
{\mathcal C}(\alpha_1 \cdot \alpha_2) &=& \{\beta_3,\beta_4\}, & X_{\beta_3} = C_1, 
X_{\beta_4} = C_2. \\  
\end{array}
$$

\item Level $2$:
The elements of level $2$ are 
$\beta_1,\beta_2, \beta_3,\beta_4,\beta_3\cdot\beta_4$. We also have,
$$
\begin{array}{ccccc}
{\mathcal C}(\beta_i) &=& \{\gamma_i\}, & X_{\gamma_i} = H_i, & i = 1,2,\\  
{\mathcal C}(\beta_i) &=& \{\gamma_i\}, & X_{\gamma_i} = C_{i-2},   &  i = 3,4,\\
{\mathcal C}(\beta_{3}\cdot\beta_{4}) &=& \{\gamma_5,\gamma_6\},& 
X_{\gamma_i} = P_{i-4},  & i = 5,6. \\  
\end{array}
$$
\end{enumerate}

We now display diagramatically the various complexes,
$\M^{\bullet,\bullet}(\alpha)$ for $\alpha \in \A_S$ starting 
at level $2$.

\begin{enumerate}
\item Level $2$:
For $1 \leq i \leq 4$, we have
\begin{eqnarray*}
\M^{\bullet,\bullet}(\beta_i)  &=&
\begin{diagram}
0 & \rTo& 0 &\rTo& 0 &\\
\uTo &&\uTo &&\uTo &\\
0 & \rTo& 0 &\rTo& 0 &\\
\uTo &&\uTo &&\uTo &\\
\HH^0(X_{\gamma_i}) & \rTo & 0&\rTo & 0&
\end{diagram}
\end{eqnarray*}

Notice that for $1 \leq i \leq 4$, 
\[
\HH^0(\Tot^{\bullet}(\M^{\bullet,\bullet}(\beta_i))) 
\cong \HH^0(X_{\beta_i}) \cong \Q.
\]

The complex $\M^{\bullet,\bullet}(\beta_3\cdot\beta_4)$ is 
shown below.

\begin{diagram}
0 & \rTo& 0 &\rTo& 0 &\\
\uTo &&\uTo &&\uTo &\\
0 & \rTo& 0 &\rTo& 0 &\\
\uTo &&\uTo &&\uTo &\\
\HH^0(P_1)\bigoplus \HH^0(P_2) & \rTo & 0&\rTo & 0&
\end{diagram}

Notice that, 
\[
\HH^0(\Tot^{\bullet}(\M^{\bullet,\bullet}(\beta_3\cdot\beta_4))) \cong \HH^0(X_{\beta_3\cdot\beta_4}) \cong \Q\bigoplus \Q.
\]

\item Level $1$:
For $i=1,2$, the complex $\M^{\bullet,\bullet}(\alpha_i)$ is as follows.
\begin{diagram}
0&\rTo&0&\rTo&0&\rTo&0& \\
\uTo&&\uTo&&\uTo&&\uTo&\\
0&\rTo&0&\rTo&0&\rTo&0& \\
\uTo&&\uTo&&\uTo&&\uTo&\\
0 & \rTo& 0 &\rTo&0&\rTo&0&\\
\uTo &&\uTo &&\uTo&&\uTo&\\
\HH^0(H_i) & \rTo & 0& \rTo & 0&\rTo & 0&
\end{diagram}

Notice that for $i=1,2$ and $j =0,1$, 
\[
\HH^j(\Tot^{\bullet}(\M^{\bullet,\bullet}(\alpha_i))) \cong 
\HH^j(H_i).
\] 

The complex $\M^{\bullet,\bullet}(\alpha_1\cdot\alpha_2)$is shown below.
\begin{diagram}
0&\rTo&0&\rTo&0&\rTo&0& \\
\uTo&&\uTo&&\uTo&&\uTo&\\
0&\rTo&0&\rTo&0&\rTo&0& \\
\uTo&&\uTo&&\uTo&&\uTo&\\
0 & \rTo& 0 &\rTo&0&\rTo&0&\\
\uTo &&\uTo &&\uTo&&\uTo&\\
\HH^0(C_1)\bigoplus\HH^0(C_2) & \rTo & \HH^0(P_1)\bigoplus\HH^0(P_2)& \rTo & 0&\rTo & 0&
\end{diagram}
Notice that for $j =0,1$, 
\[
\HH^j(\Tot^{\bullet}(\M^{\bullet,\bullet}(\alpha_1\cdot\alpha_2))) \cong  \HH^j(H_{12}).
\]

\item Level $0$:

The complex $\M^{\bullet,\bullet}(0)$ is shown below:

\begin{diagram}
0&\rTo&0&\rTo&0&\rTo&0&\rTo&0& \\
\uTo&&\uTo&&\uTo&&\uTo&&\uTo&\\
0&\rTo&0&\rTo&0&\rTo&0&\rTo&0& \\
\uTo&&\uTo&&\uTo&&\uTo&&\uTo&\\
0&\rTo&0&\rTo&0&\rTo&0&\rTo&0& \\
\uTo&&\uTo&&\uTo&&\uTo&&\uTo&\\
0 & \rTo& \HH^0(P_1)\bigoplus\HH^0(P_2) &\rTo&0&\rTo&0&\rTo&0&\\
\uTo &&\uTo^{d^{1,0}} &&\uTo&&\uTo&&\uTo&\\
\HH^0(H_1)\bigoplus\HH^0(H_2) & \rTo^{\delta^{0,0}} & \HH^0(C_1)\bigoplus\HH^0(C_2)& \rTo & 0&\rTo & 0&\rTo & 0&
\end{diagram}

The matrices for the homomorphisms, $\delta^{0,0}$ and $d^{1,0}$ in the
obvious bases are both equal to
$$
\left(
\begin{array}{cc}
1\;\; &\;\; 1 \\
1\;\; &\;\; 1
\end{array}
\right).
$$

From the fact that the rank of the above matrix is $1$, 
it is not too difficult to deduce that, 
$\HH^j(\Tot^{\bullet}(\M^{\bullet,\bullet}(0))) \cong \HH^j(S)$,
for $j =0,1,2$, that is
$$
\begin{array}{ccc}
\HH^0(\Tot^{\bullet}(\M^{\bullet,\bullet}(0))) &\cong & \Q, \cr
\HH^1(\Tot^{\bullet}(\M^{\bullet,\bullet}(0))) &\cong & 0, \cr
\HH^2(\Tot^{\bullet}(\M^{\bullet,\bullet}(0))) &\cong &  \Q. 
\end{array}
$$
\end{enumerate}
\end{example}

We now prove properties (1) and (2) of the various 
$\M^{\bullet,\bullet}(\alpha)$.

\begin{proposition}
\label{prop:main}
For each $\alpha\in \A_S$ 
the double complex $\M^{\bullet,\bullet}(\alpha)$ satisfies the following 
properties:

\begin{enumerate}
\item
$
\HH^i(\Tot^{\bullet}(\M^{\bullet,\bullet}(\alpha)))  \cong 
\HH^i(X_\alpha)$ 
for $0 \leq i \leq \ell - \level(\alpha)$.
\item
For every $\beta \in \A_S$, such that $\alpha$ is an ancestor of $\beta$,
and $\level(\alpha) = \level(\beta)$,
the homomorphism,
$ r_{\alpha,\beta}^{\bullet,\bullet}:  
\M^{\bullet,\bullet}(\alpha) \rightarrow \M^{\bullet,\bullet}(\beta),$
induces the  restriction homomorphisms between the cohomology
groups:
$$
\displaylines{
r^*: \HH^i(X_\alpha) \longrightarrow \HH^i(X_\beta)
}
$$ 
for $0 \leq i \leq \ell - \level(\alpha)$  via the isomorphisms in (1).
\end{enumerate}
\end{proposition}

The main idea behind the proof of Proposition \ref{prop:main} 
is as follows.
We consider a triangulation $h_0: \Delta_0 \rightarrow S$, such that
for any $\alpha \in A_S$, 
$h_0$ restricts to a semi-algebraic triangulation,
$h_\alpha: \Delta_\alpha \rightarrow X_\alpha$.
Note that, this implies that if
$\beta \in \A_S$ and  $\alpha \in \ancestor(\beta)$,
then the triangulation
$h_\alpha: \Delta_{\alpha} \rightarrow X_\alpha$ 
restricts to the 
triangulation  $h_\beta: \Delta_{\beta} \rightarrow X_{\beta}$,
and in particular $\Delta_\beta$ is a subcomplex of $\Delta_\alpha$.


For each $\alpha \in \A_S$, we have that 
$\Delta_\alpha = \cup_{\alpha_0 \in {\mathcal C}(\alpha)} \Delta_{\alpha_0}$,
and each $\Delta_{\alpha_0}$ for $\alpha_0 \in {\mathcal C}(\alpha)$ is
a subcomplex of $\Delta_\alpha$.
We denote by $\N^{\bullet,\bullet}(\Delta_{\alpha})$ the Mayer-Vietoris
double complex of $\Delta_{\alpha}$ with respect to the sub-complexes 
$\Delta_{\alpha_0}, \alpha_0 \in {\mathcal C}(\alpha)$ (cf. Definition
\ref{def:MV}).

We denote by $n_\alpha = \ell - \level(\alpha) +1$.
Recall that $\N_{n_\alpha}^{\bullet,\bullet}(\Delta_{\alpha})$
is the following truncated complex.

\[
\begin{array}{ccll}
\N_{n_\alpha}^{p,q}(\Delta_{\alpha}) & = & \N^{p,q}(\Delta_{\alpha}),\; & 
\;0 \leq p+q \leq n_\alpha, \cr
\N_{n_\alpha}^{p,q}(\Delta_{\alpha}) & = &  0,\; &\; \mbox{otherwise}. \cr
\end{array}
\]

By Corollary \ref{cor:MV} we have that,
$$
\displaylines{
\HH^i(\Tot^{\bullet}(\N_{n_\alpha}^{\bullet,\bullet}(\Delta_{{\alpha}}))) 
\cong \HH^i(X_{\alpha}), \; 0 \leq i \leq \ell - \level(\alpha).
}
$$

We then prove by induction on $\level(\alpha)$ that for each 
$\alpha \in \A_S$  there exists a double complex $D^{\bullet,\bullet}(\alpha)$ and 
homomorphisms,
$$
\displaylines{
\phi_\alpha^{\bullet,\bullet}: \M^{\bullet,\bullet}(\alpha) 
\longrightarrow D^{\bullet,\bullet}(\alpha) \cr
\psi_\alpha^{\bullet}: \Ch^{\bullet}(\Delta_{\alpha}) \longrightarrow \Tot^{\bullet}
(D^{\bullet,\bullet}(\alpha)) 
}
$$
such that,
$$
\displaylines{
\Tot^{\bullet}(\phi_\alpha^{\bullet,\bullet}): \Tot^{\bullet}(\M^{\bullet,\bullet}(\alpha))
 \longrightarrow \Tot^{\bullet}(D^{\bullet,\bullet}(\alpha)), 
}
$$
as well as $\psi_\alpha^{\bullet}$ (as shown in the following figure)
are quasi-isomorphisms.

\begin{diagram}
&& \Tot^{\bullet}(D^{\bullet,\bullet}(\alpha)) && \cr
&\ruTo^{\Tot^{\bullet}(\phi_\alpha^{\bullet,\bullet})} && 
\luTo^{\psi_\alpha^{\bullet}}& \cr
\Tot^{\bullet}(\M^{\bullet,\bullet}(\alpha)) &&&&
\Ch^{\bullet}(\Delta_{\alpha})
\end{diagram}

These quasi-isomorphisms will together imply that,
\[
\HH^i(\Tot^{\bullet}(\M^{\bullet,\bullet}(\alpha))) \cong 
\HH^i(\Tot^{\bullet}(D^{\bullet,\bullet}(\alpha))) \cong 
\HH^i(\Tot^{\bullet}(\N_{n_\alpha}^{\bullet,\bullet}(\Delta_{\alpha}))) \cong 
\HH^i(X),
\]
for $0 \leq i \leq \ell - \level(\alpha).$

\begin{proof}{Proof of Proposition \ref{prop:main}:}
The proof of the proposition is by induction on $\level(\alpha)$. 
When $\level(\alpha) = \ell$, we let
$D^{\bullet,\bullet}(\alpha) = 
\N_{n_\alpha}^{\bullet,\bullet}(\Delta_{\alpha})$, and define the
homomorphisms $\phi_\alpha^{\bullet,\bullet},
\psi_\alpha^{\bullet}$ as follows.
From the definition of 
$\M^{\bullet,\bullet}(\alpha)$ it is clear that in order to define
$\phi_\alpha^{\bullet,\bullet}$, it suffices to define,
$\phi_\alpha^{0,0}$ and $\phi_\alpha^{0,1}$.

We define,
\[
\phi_\alpha^{0,0}: \M^{0,0}(\alpha) = 
\bigoplus_{{\alpha_0} \;\in\; {\mathcal C}(\alpha)} 
\;\HH^0(X_{\alpha_0})
\rightarrow
\bigoplus_{{\alpha_0} \;\in\; {\mathcal C}(\alpha)} \Ch^0(\Delta_{{\alpha_0}})
= \N_1^{0,0}(\Delta_{X_\alpha}),
\]
by defining for 
$\theta \in 
\bigoplus_{{\alpha_0} \;\in\; {\mathcal C}(\alpha)} 
\;\HH^0(X_{\alpha_0})$, and
any vertex $v$ of the complex $\Delta_{{\alpha_0}}$,
$\phi_\alpha^{0,0}(\theta)_{\alpha_0}(v)$ to be the value of the 
locally constant function  $\theta_{\alpha_0}$ on 
$X_{\alpha_0}$.

Similarly, we define 
\[
\phi_\alpha^{0,1}: 
\bigoplus_{{\alpha_0 < \alpha_1,}\alpha_i \in {\mathcal C}(\alpha)} 
\;\HH^0(X_{\alpha_0\cdot\alpha_1})
\rightarrow
\bigoplus_{{\alpha_0 <\alpha_1} \alpha_i\in{\mathcal C}(\alpha)} \Ch^0(\Delta_{{\alpha_0\cdot\alpha_1}}),
\]
noting that
\[
\M^{0,1}(\alpha) = \bigoplus_{\alpha_0 < \alpha_1,\alpha_i\in {\mathcal C}(\alpha)} 
\;\HH^0(X_{\alpha_0\cdot\alpha_1}),
\]
and
\[
\N_1^{0,0}(\Delta_{\alpha}) = 
\bigoplus_{\alpha_0 <\alpha_1, \alpha_i \in {\mathcal C}(\alpha)} \Ch^0(\Delta_{{\alpha_0\cdot\alpha_1}}),
\]
by defining for
$\theta \in 
\bigoplus_{\alpha_0 < \alpha_1,\alpha_i \in {\mathcal C}(\alpha)} 
\;\HH^0(X_{\alpha_0\cdot\alpha_1})
$, and
any vertex $v$ of the complex $\Delta_{{\alpha_0\cdot\alpha_1}}$,
$\phi_\alpha^{0,1}(\theta)_{\alpha_0,\alpha_1}(v)$ to be the value of the 
locally constant function  $\theta_{\alpha_0,\alpha_1}$ on 
the connected component of $X_{\alpha_0\cdot\alpha_1}$ containing
$h_{\alpha_0\cdot\alpha_1}(v)$.

The homomorphism $\psi^{\bullet}_{\alpha}$ is induced by restriction
as in the definition of $\psi^{\bullet}_{\ell+1}$ 
in Corollary \ref{cor:MV}.

It is now easy to verify, that 
$\Tot^{\bullet}(\phi_\alpha^{\bullet,\bullet})$ and
$\psi_\alpha^{\bullet}$
are indeed quasi-ismorphisms.

In general for $\alpha \in \A_S$, with $\level(\alpha) < \ell$,
we have by induction that for
each ${\alpha_0},\ldots,{\alpha_p}, {\alpha_{p+1}} \in {\mathcal C}(\alpha),
0 \leq p \leq \ell - \level(\alpha) + 2,$
there exists a double complex 
$D^{\bullet,\bullet}({\alpha_0\cdots\alpha_p})$ and 
quasi-isomorphisms
$$
\displaylines{
\Tot^{\bullet}(\phi_{{\alpha_0\cdots\alpha_p}}^{\bullet,\bullet}): 
\Tot^{\bullet}(\M^{\bullet,\bullet}({\alpha_0\cdots\alpha_p}))
 \longrightarrow \Tot^{\bullet}(D^{\bullet,\bullet}({\alpha_0\cdots\alpha_p}))
 \cr
\psi_{{\alpha_0\cdots\alpha_p}}^{\bullet}: 
\Ch_{n_\alpha}^{\bullet}(\Delta_{{\alpha}})
\longrightarrow 
\Tot^{\bullet}(D^{\bullet,\bullet}({\alpha_0\cdots\alpha_p})) .
}
$$

We now define $D^{\bullet,\bullet}(\alpha)$ by,
\[
\begin{array}{cccc}
D^{p,q}(\alpha) &=&   \bigoplus_{\alpha_0< \cdots <\alpha_p, \;
{\alpha_i} \in {\mathcal C}(\alpha)}
\; 
\Tot^q(D^{\bullet,\bullet}({\alpha_0\cdots\alpha_p})), \; & \; 
0 \leq  p+q \leq n_{\alpha}, \\
&=& 0, \;&\; \mbox{else}.
\end{array}
\]
The homomorphism $\phi_\alpha^{\bullet,\bullet}$ is the one induced by the
different $\Tot^{\bullet}(\phi_{{\alpha_0\cdots\alpha_p}}^{\bullet,\bullet})$ 
defined already by induction, that is
\[
\phi_{\alpha}^{p,q}: \M^{p,q}(\alpha) \rightarrow D^{p,q}(\alpha),
\]
is defined by
\[
\phi_{\alpha}^{p,q} = \bigoplus_{\alpha_0< \cdots <\alpha_p, \;
{\alpha_i} \in {\mathcal C}(\alpha)}\Tot^q(\phi^{\bullet,\bullet}_{\alpha_0\cdots\alpha_p}).
\]

In order to define the  homomorphism $\psi_\alpha^{\bullet}$,
we first define a homomorphism,
\[
\rho_\alpha^{\bullet,\bullet}: \N_{n_\alpha}^{\bullet,\bullet}(\Delta_{{\alpha}}) 
\longrightarrow D^{\bullet,\bullet}(\alpha)
\]
induced by the  different $\psi_{{\alpha_0\cdots\alpha_p}}^{\bullet}.$

We define
\[
\rho_{\alpha}^{p,q}: \N_{n_\alpha}^{p,q}(\Delta_{\alpha}) \rightarrow D^{p,q}(\alpha),
\]
by
\[
\rho_{\alpha}^{p,q} = \bigoplus_{\alpha_0< \cdots <\alpha_p, \;
{\alpha_i} \in {\mathcal C}(\alpha)}\psi_{\alpha_0\cdots\alpha_p}^{q}.
\]

We now compose the homomorphism,
\[
\Tot^{\bullet}(\rho_\alpha^{\bullet,\bullet}): 
\Tot^{\bullet}(\N_{n_\alpha}^{\bullet,\bullet}(\Delta_{{\alpha}}) )
\longrightarrow 
\Tot^{\bullet}(D^{\bullet,\bullet}(\alpha)),
\] 
with the quasi-isomorphism 
\[
\psi_{\alpha, n_{\alpha}}^{\bullet}: \Ch_{n_\alpha}^{\bullet}(\Delta_{{\alpha}}) \longrightarrow 
\Tot^{\bullet}(\N_{n_\alpha}^{\bullet,\bullet}(\Delta_{{\alpha}}))
\]
(see Proposition \ref{prop:MV}). 

Using the induction hypothesis it is easy to see that the
homomorphism $\phi_{\alpha}^{\bullet,\bullet}$ induces an isomorphism
between the $'E_1$ terms of the corresponding spectral sequences. It follows
from Proposition \ref{prop:spectral} that this implies that $\Tot^{\bullet}(\phi_\alpha^{\bullet,\bullet})$ is a quasi-ismorphism.
A similar argument also shows that $\Tot^{\bullet}(\rho_\alpha^{\bullet,\bullet})$ is also a quasi-isomorphism and hence so is $\psi_\alpha^{\bullet}$
since it is a composition of two quasi-isomorphisms.
This completes the induction.
\end{proof}

\section{Algorithmic Preliminaries}
\label{sec:algo_prelim}
In this section, we describe some algorithmic results which we need in
the main algorithms.

\subsection{Computation with Complexes}
\label{subsec:linearalgebra}
In the description of our algorithm, we compute in a recursive way
certain complicated double complexes, whose constructions have already
been described in Section \ref{sec:covering}. 
The computation of a complex (or a
double complex) means computing bases for each term of the complex
(or double complex), as well as the matrices representing the differentials
in this bases. Given a complex $\Ch^{\bullet}$ (in terms of some fixed bases),
we can compute its homology groups $\HH^*(\Ch^{\bullet})$ using elementary
algorithms from linear algebra for computing kernels and images of 
vector space homomorphisms. Similarly, given a double complex,
$\D^{\bullet,\bullet},$ we can compute the 
complex $\Tot^{\bullet}(D^{\bullet,\bullet})$ as well as,
$\HH^*(\Tot^{\bullet}(D^{\bullet,\bullet}))$, using linear algebraic 
subroutines. Since the naive algorithms (using say Gaussian elimination
for computing kernels and images of linear maps) run in time polynomial
in the dimensions of the vector spaces involved, it is clear that all the
above computations involving complexes can be done in time polynomial
in the sum of the dimensions of all terms in the input complex. 
This is sufficient for proving the main result of this paper,
and we do not make any attempt to perform these computations in an
optimal manner using more sophisticated algorithms.

\subsection{General Position and Covers by Contractible Sets}
\label{subsec:covers}
We first recall some results  proved in \cite{BPR04} on constructing
singly exponential sized cover of a given closed semi-algebraic set,
by closed, contractible semi-algebraic set.
We recall the input, output and the complexity of the algorithms, referring
the reader to \cite{BPR04} for all details including the proofs of
correctness.

\subsubsection{General Position}
Let $Q \in \R[X_1,\ldots,X_k]$
such that
$\ZZ(Q,\R^k) = \{ x \in \R^k \mid Q(x) = 0 \}$ is bounded.
We say that a finite set of
 polynomials ${\mathcal P} \subset \D[X_1,\ldots,X_k]$
is in strong $\ell$-general position with respect to $Q$
if
any  $\ell+1$
polynomials belonging to 
${\mathcal P}$ have no  zeros in common with
$Q$ in $\R^k$, and 
any  $\ell$
polynomials belonging to 
${\mathcal P}$ have at  most a finite number of  zeros in common with
$Q$ in $\R^k$. 

\subsubsection{Infinitesimals}
In our algorithms we will use infinitesimal perturbations.
In order to do so, we will extend the ground field $\R$ 
to,
$\R\langle \varepsilon\rangle$,  the real closed field of algebraic
Puiseux series in $\varepsilon$ with coefficients in $\R$ \cite{BPR03}. 
The sign of a Puiseux series in $\R\langle \varepsilon\rangle$
agrees with the sign of the coefficient
of the lowest degree term in
$\varepsilon$. 
This induces a unique order on $\R\langle \varepsilon\rangle$ which
makes $\varepsilon$
infinitesimal: $\varepsilon$ is positive and smaller than
any positive element of $\R$.
When $a \in \R\la \varepsilon \ra$ is bounded by an element of $\R$,
$\lim_\varepsilon(a)$ is the constant term of $a$, obtained by
substituting 0 for $\varepsilon$ in $a$.
We will also denote the field 
$\R\la\varepsilon_1\ra\cdots\la\varepsilon_s\ra$ by $\R\la\bar\varepsilon\ra$,
where $\varepsilon_1, \varepsilon_2, \ldots,\varepsilon_s>0$ are 
infinitesimals with respect to the field $\R$.

\subsubsection{Replacement by closed sets without changing cohomology}
The following  algorithm allows us to replace a given semi-algebraic set by
a new one which is closed and defined by polynomials in general position
and which has the same homotopy type as the the given set. 
This construction is essentially due
to Gabrielov and Vorobjov \cite{GV}, where it was shown that the
sum of the Betti numbers is preserved. The homotopy equivalence property is
shown in \cite{BPR04}.

\begin{algorithm2} [Cohomology Preserving Modification to Closed]
\label{alg:closed}
\begin{description}
\item[]
\item [{\bf Input} :]
\begin{enumerate}
\item an element $c \in \R$, such that $c >0$,
\item
a polynomial $Q \in \R[X_1,\ldots,X_k]$ such that $\ZZ(Q,\R^k)\subset 
B(0,1/c)$, 
\item  
a finite set of $s$ polynomials 
$$
{\mathcal P} =
\{P_1,\ldots,P_s \} \subset \R[X_1,\ldots,X_k],
$$
\item
a subset $\Sigma\subset {\rm Sign}(Q,{\mathcal P})$, defining
a semi-algebraic set $X$ by
$$
\displaylines{
X = \cup_{\sigma \in \Sigma}\RR(\sigma).
}
$$
\end{enumerate}
\item [{\bf Output} :]
A description of a 
${\mathcal P}'$-closed and bounded semi-algebraic subset, 
$$
X' \subset \ZZ(Q,\R\langle\eps, \varepsilon_1,\ldots,\varepsilon_{2s}\rangle^k),
$$ 
with
${\mathcal P}' = \bigcup_{1 \leq i \leq s, 1 \leq j \leq 2s} 
\{P_i \pm \varepsilon_j\},$
such that,
\begin{enumerate}
\item
$\HH^*(X')\cong \HH^*(X)$, and
\item
the family of polynomials ${\mathcal P}'$ is in $k'$-strong general position
with respect to $\ZZ(Q,\R\langle\eps,\varepsilon_1,\ldots,\varepsilon_{2s}
\rangle^k)$,
where $k'$ is the real dimension of 
\[
\ZZ(Q,\R\langle\eps,\varepsilon_1,\ldots,\varepsilon_{2s}\rangle^k).
\]
\end{enumerate}
\item [{\bf Procedure} :]
\item[Step 1] Let $\eps$ be an infinitesimal. 
\begin{enumerate}
\item
Define $\tilde T$ as the intersection of
$\E(T,\R\la\eps\ra)$ with the ball of center $0$ and radius $1/\eps$.
\item
Define ${\mathcal P}$ as ${\mathcal Q} \cup \{\eps^2(X_1^2+\ldots+X_k^2+X_{k+1}^2)-4, X_{k+1}\}$. 
\item 
Replace $\tilde T$ by the ${\mathcal P}$- semi-algebraic set $S$ defined as the intersection
 of the cylinder $\tilde T\times  \R\la \eps \ra$ with the upper hemisphere defined by
 $\eps^2(X_1^2+\ldots+X_k^2+X_{k+1}^2)=4, X_{k+1}\ge 0$.
\end{enumerate}
\item[Step 2]
Using the Gabrielov-Vorobjov construction described in \cite{BPR04},   
replace $S$ by   a ${\mathcal P}'$-closed set, $S'$.
Note that ${\mathcal P}'$ is in general
position with respect to the sphere of center $0$ and radius $2/\eps$.
\end{description}
\end{algorithm2}

\begin{proof}{Complexity:}
Let $d$ be the maximum degree among the polynomials in ${\mathcal P}$.
The total complexity is bounded by $s^{k+1}d^{O(k)}$ (see \cite{BPR04}).
\end{proof}

\subsubsection{Algorithm for Computing Covers by Contractible Sets}
The following algorithm described in detail in \cite{BPR04} is used to 
a cover of a given closed and bounded
semi-algebraic sets defined by polynomials in general position by
closed, bounded and contractible semi-algebraic sets.

\begin{algorithm2} [Cover by Contractible Sets]
\label{16:alg:acycliccovering}
\index{Connecting!Algorithm!Acyclic Covering}
\begin{description}
\item[]
\item [{\bf Input} :] 
\begin{enumerate}
\item 
a polynomial $Q \in \D[X_1,\ldots,X_k]$ such that $\ZZ(Q,\R^k)\subset 
B(0,1/c)$,
\item  
a finite set of $s$ polynomials ${\mathcal P} \subset \D[X_1,\ldots,X_k]$
in strong $\ell$-general position on $\ZZ(Q,\R^k)$.
\end{enumerate}
\item [{\bf Output} :]
\begin{enumerate}
\item
a finite  family of polynomials ${\mathcal C}
=  \{Q_1,\ldots,Q_N \}
\subset \R[X_1,\ldots,X_k],$ 
\item
the finite family $\overline{\mathcal C} \subset 
\R[\bar\varepsilon][X_1,\ldots,X_k]$ 
(where $\bar\eps$ denotes the infinitesimals
$\varepsilon_1 \gg \varepsilon_2 \gg \cdots \gg \varepsilon_{2N} > 0$)
defined by
$$\overline{\mathcal C}= \{Q\pm \varepsilon_i \mid Q \in {\mathcal C}, 1\le i \le 2N\}.$$
\item
a set of $\overline{\mathcal C}$-closed formulas $\{\phi_1,\ldots,\phi_M\}$ 
such that
\begin{enumerate}
\item 
each  $\RR(\phi_i,\R\langle\bar\varepsilon\rangle^k)$ is contractible, 
\item 
their union
$\displaystyle{
\cup_{1 \leq i \leq M} 
\RR(\phi_i,\R\langle\bar\varepsilon\rangle^k) = 
\ZZ(Q,\R\langle\bar\varepsilon\rangle^k), and 
}
$
\item
each basic ${\mathcal P}$-closed subset of 
$\ZZ(Q,\R\langle\bar\varepsilon\rangle^k)$ is a union of some subset of the
$\RR(\phi_i,\R\langle\bar\varepsilon\rangle^k)$'s.
\end{enumerate}
\end{enumerate}

\end{description}
\end{algorithm2}

\begin{proof}{Complexity:}
The total complexity is bounded by $s^{(k+1)^2}d^{O(k^5)}$ (see \cite{BPR04}).
\end{proof}

\section{Algorithm for computing the first $\ell$ Betti numbers of a 
semi-algebraic set}
\label{sec:main}
We are finally in a position to describe the main algorithm of this paper.
\begin{algorithm2}
\label{algo:main}
[First $\ell$ Betti Numbers of a ${\mathcal P}$ Semi-algebraic Set]
\begin{description}
\item[]
\item [{\bf Input} :] 
\begin{description}
\item 
a polynomial $Q \in \D[X_1,\ldots,X_k]$ such that $\ZZ(Q,\R^k)\subset B(0,1/c)$,
\item  
a finite set of  polynomials ${\mathcal P} \subset \D[X_1,\ldots,X_k],$
\item a formula defining a ${\mathcal P}$ semi-algebraic set $S$ contained in
$\ZZ(Q,\R^k)$.
\end{description}

\item [{\bf Output} :]
$b_0(S),\ldots,b_\ell(S).$
\item [{\bf Procedure} :]
\item[Step 1]
Using Algorithm \ref{alg:closed} (Cohomology Preserving Modification 
to Closed),   
replace $S$ by   a ${\mathcal P}'$-closed set, $S'$.
Note that ${\mathcal P}'$ is in $k'$-general position 
with respect to $\ZZ(Q,\R^k)$.

\item[Step 2]
Use Definition \ref{def:admissible} to compute $\A_{S'}$ using 
Algorithm \ref{16:alg:acycliccovering} (Cover by Contractible Sets)
for computing the various ${\mathcal C}'(\cdot)$ occuring in the
definition of $\A_{S'}$. 
For each element $\alpha \in \A_{S'}$, we also compute the set
of ancestors $\ancestor(\alpha) \subset \A_{S'}$,
${\mathcal C}(\alpha)$, as well as  $\level(\alpha)$.

More precisely, we do the following.
\begin{enumerate}
\item
\begin{enumerate}
\item
Initialize,
$$
\displaylines{
\A_{S'} \leftarrow \emptyset,
}
$$
\item
$$
\displaylines{
\A_{S'} \leftarrow  \A_{S'} \cup \{0\}, \cr
\level(0)\leftarrow  0, \cr
X_0 \leftarrow  S', \cr
{\mathcal C}(0) \leftarrow {\mathcal C}'(S'), \cr
\ancestor(0) = \{0\}.
}
$$
Also, maintain a directed graph $G$ with the current set $\A_{S'}$ as
its set of vertices representing the ancestor-descendent relationships.
\end{enumerate}

\item
For $i=0$ to $\ell$ do the following:
\begin{enumerate}
\item
For each $\alpha \in \A_{S'}$ at level $i$, 
with $\ancestor(\alpha) = \{\alpha_1,\ldots,\alpha_N\}$,

$$
\displaylines{
{\mathcal C}(\alpha) \leftarrow  \bigcup_{\beta_i \in {\mathcal C}(\alpha_i), 1 \leq i \leq N}
{\mathcal C}'(\beta_1  \cdots \beta_N \cdot \alpha)
}
$$
using Algorithm \ref{16:alg:acycliccovering} (Cover by Contractible Sets).
\item
For $0 \leq j \leq \ell-i+1$ and
each $\alpha_0,\ldots,\alpha_j  \in {\mathcal C}(\alpha)$,
\[
A_{S'} \leftarrow A_{S'} \cup \{\alpha_0\cdot\alpha_1\cdots\alpha_j\},
\]

\[
X_{\alpha_0\cdots\alpha_j} \leftarrow  
X_{\alpha_0} \cap \cdots \cap X_{\alpha_j},
\]  

\[
\level(\alpha_0\cdot\alpha_1\cdots\alpha_j) \leftarrow i+1.
\]
\item
For each 
$\{\alpha_0,\ldots,\alpha_i\} \subset \{\beta_0,
\ldots,\beta_j \} \subset {\mathcal C}(\alpha)$,  
with $j \leq \ell-i+1$, 
\[
 \ancestor(\beta_0\cdots\beta_j) \leftarrow
\ancestor(\beta_0\cdots\beta_j) \cup \{\alpha_0\cdots\alpha_i\},
\] 
and update $G$.
\item
For each 
$\alpha' \in {\mathcal C}'(\beta_1 \cdot \cdots\cdot \beta_N 
\cdot \alpha)$, 
\[\ancestor(\alpha') \leftarrow
\ancestor(\alpha') \cup  \{\beta_1,\ldots,\beta_N\}.
\]
and update $G$.
Use any graph transitive closure algorithm to transitively close $G$.
Accordingly  update all the sets $\ancestor(\alpha), \alpha \in \A_{S'}$. 
\end{enumerate}
\end{enumerate}
\item[Step 3]
Using Definition \ref{def:double}, compute for each $\alpha \in \A_{S'}$,
the complex $\M^{\bullet,\bullet}(\alpha)$ starting with elements
$\alpha \in \A_{S'}$ with $\level(\alpha) = \ell$. 
Note that for each $\alpha \in \A_{S'}$, ${\mathcal C}(\alpha)$ has 
already being computed in Step 2.
This allows us to compute 
matrices corresponding to all the homomorphisms in
$\M^{\bullet,\bullet}(\alpha)$ for $\alpha \in \A_{S'}$ with
$\level(\alpha) = \ell.$
The recursive definition
of $\M^{\bullet,\bullet}(\alpha)$, implies that we can compute the matrices 
corresponding to all the homomorphisms in
$\M^{\bullet,\bullet}(\alpha)$ for 
$\alpha \in \A_{S'}$ with
$\level(\alpha) < \ell$,
once we have computed the same for $\M^{\bullet,\bullet}(\beta)$, for 
all $\beta \in \A_{S'}$ with $\level(\beta) > \level(\alpha).$
The same is also true for the matrices corresponding to the restriction
homomorphisms $r_{\alpha,\beta}^{\bullet,\bullet}$.
\item[Step 4]
For each $i, 0 \leq i \leq \ell,$ compute
\[
b_i(S) = \dim_{\Q} \HH^i(\Tot^{\bullet}({\mathcal M}^{\bullet,\bullet}(0))),
\]
using standard linear algebra algorithms for computing dimensions of 
kernels and images of linear transformations. 
\end{description}
\end{algorithm2}

\begin{proof}{Proof of correctness :}
The correctness of the algorithm is a consequence of the correctness of
Algorithms \ref{alg:closed} (Cohomology Preserving Modification to 
Closed), 
Algorithm \ref{16:alg:acycliccovering} (Cover by Contractible Sets),
and Proposition \ref{prop:main}.
\end{proof}

\begin{proof}{Complexity analysis:}
\hide{
Each step is clearly singly exponential from the complexity analysis
of Algorithms \ref{alg:closed} (Cohomology Preserving Modification to 
Closed), \ref{16:alg:acycliccovering} (Cover by Contractible Sets), and
Proposition \ref{prop:bound}.
}
The complexity of Step 1 is bounded by $(sd)^{O(k)}$ using the complexity
analysis of Algorithm \ref{alg:closed} (Cohomology Preserving Modification 
to Closed).

In order to bound the complexity of Step 2, note that the number
of calls to 
Algorithm \ref{16:alg:acycliccovering} (Cover by Contractible Sets).
for computing various covers,
${\mathcal C}'(\cdot)$ is bounded by $\#\A_{S'}$,which in turn is 
bounded by $(sd){k^{O(\ell)}}$ by Proposition \ref{prop:bound}. Moreover,
the cost of each such call is also bounded by $(sd){k^{O(\ell)}}$. The cost of
all other operations, including updating the list of ancestors of elements of
$\A_{S'}$ is polynomial in $\#\A_{S'}$. Thus, the total complexity of this 
step is bounded by $(sd){k^{O(\ell)}}$.
Finally,
the complexity of the computations involving linear algebra
in Step 3 is polynomial in the cost of
computing the various complexes $\M^{\bullet,\bullet}(\alpha)$, as well their
sizes (see Section \ref{subsec:linearalgebra}). 
All these are bounded by $(sd)^{k^{O(\ell)}}$ by Proposition
\ref{prop:bound}.
Thus, the complexity of the whole algorithm is bounded by
$(sd)^{k^{O(\ell)}}$.
\end{proof}

\section{Implementation and Practical Aspects}
\label{sec:practical}
The problem of computing all the Betti numbers of semi-algebraic sets in 
single exponential time (as well as the related problems of existence
of single exponential sized triangulations or even stratifications)
is considered a very important
question in quantitative real algebraic geometry. The main result of this 
paper should be considered a partial progress on this theoretical problem. 
Since
the complexity of Algorithm \ref{16:alg:acycliccovering} 
(Cover by Contractible Sets) for computing contractible covers
is very high (even though single exponential), the complexity of Algorithm
\ref{algo:main} is prohibitively expensive for practical implementation.
The topological ideas underlying our algorithm has been implemented in 
a very limited setting in order to compute the first two Betti numbers
of sets defined by quadratic inequalities (see \cite{BK05}).  
In this implementation,
the covering is obtained by means different from 
Algorithm \ref{16:alg:acycliccovering}. 
However, practical implementation for general
semi-algebraic sets remains a formidable challenge.


\begin{thebibliography}{50}
\bibitem[{B99}]{B99}
{\sc S.\ Basu},
\newblock {\em On Bounding the Betti Numbers and Computing the Euler
Characteristics of Semi-algebraic Sets},
\newblock {Discrete and Computational Geometry}, 22 1-18 (1999).

\bibitem[{B03}]{B03}
{\sc S.\ Basu},
\newblock{\em On different bounds on different Betti numbers},
\newblock {Discrete and Computational Geometry}, Vol 30, No. 1 (2003).

\bibitem[{BK05}]{BK05}
{\sc S.\ Basu, M.\ Kettner},
\newblock{ \em Computing the Betti numbers of arrangements in practice,}
\newblock Proceedings of the 8-th International Workshop
on Computer Algebra in Scientific Computing (CASC), LNCS 3718, 13-31, 2005.


\bibitem[{BPR95}]{BPR95}
{\sc S.\ Basu, R.\ Pollack, M.-F. \ Roy},
 \newblock {\em On the Combinatorial and Algebraic 
Complexity of Quantifier Elimination},
\newblock Journal of the ACM , 43  1002--1045, (1996).

\bibitem[{BPR99}]{BPR99}
{\sc S.\ Basu, R.\ Pollack, M.-F. \ Roy},
 \newblock {\em Computing Roadmaps of Semi-algebraic
Sets on a Variety},
\newblock Journal of the AMS, vol 3, 1 55-82 (1999).

\bibitem[{BPR03}]{BPR03}
{\sc S.\ Basu, R.\ Pollack, M.-F. \ Roy},
\newblock {\em Algorithms in Real Algebraic Geometry},
Springer-Verlag, 2003.

\bibitem[{BPR04}]{BPR04}
{\sc S.\ Basu, R.\ Pollack, M.-F. \ Roy},
 \newblock {\em Computing the first Betti number and the connected components
of semi-algebraic sets}, preprint (2004). 
\newline 
(Available at {\tt www.math.gatech.edu/$\tilde{}$saugata/bettione.ps}.)


\bibitem[{BCR}]{BCR}
{\sc J.\ Bochnak, M.\ Coste, M.-F.\ Roy},
\newblock {\em G\'eom\'etrie
alg\'ebrique r\'eelle,}
Springer-Verlag (1987).
\newblock {\em Real algebraic geometry},
Springer-Verlag (1998).

\bibitem[{BC04}]{BC}
{\sc P.\ Burgisser, F.\ Cucker},
\newblock {\em Counting Complexity Classes for Numeric Computations II: 
Algebraic and  Semi-algebraic Sets},
\newblock preprint.

\bibitem[{Canny93}]{Canny93a}
{\sc J.\ Canny},
\newblock {\em Computing road maps in general semi-algebraic sets},
\newblock The
Computer Journal, 36: 504--514, (1993).

\bibitem[{Collins75}]{Col}
{\sc G. Collins},
\newblock {\em Quantifier elimination for real closed fields by
cylindric algebraic decomposition},
\newblock In
Second GI Conference on Automata Theory and
Formal Languages. Lecture
Notes in Computer Science, vol. 33, pp. 134-183, Springer-
Verlag, Berlin (1975).

\bibitem[{GV05}]{GV}
{\sc A.\ Gabrielov, N.\ Vorobjov}
\newblock {Betti Numbers for Quantifier-free Formulae,}
\newblock Discrete and Computational Geometry, 33:395-401, 2005.

\bibitem[{GR92}]{GR92}
{\sc L. Gournay, J. J. Risler},
\newblock{\em Construction of roadmaps of semi-algebraic sets},
\newblock
Appl. Algebra Eng. Commun. Comput. 4, No.4, 239-252 (1993).


\bibitem[{GV92}]{GV92}
{\sc D. Grigor'ev, N. Vorobjov},
\newblock{\em Counting connected components of a semi-algebraic
set in subexponential time},
\newblock Comput. Complexity 2, No.2, 133-186 (1992).

\bibitem[{Hardt80}]{Hardt} {\sc R. M. Hardt},
\newblock{\em Semi-algebraic Local Triviality in Semi-algebraic Mappings},
\newblock Am. J. Math. { 102}, 291-302 (1980).

\bibitem[{HRS94}]{HRS94}
{\sc J.\ Heintz, M.-F.\ Roy, P. Solern\`o},
\newblock{\em Description of the
Connected Components of a Semialgebraic Set in Single Exponential Time},
\newblock{Discrete and Computational Geometry}, 11, 121-140 (1994).

\bibitem[{Knebusch89}]{Knebusch}
{\sc M.\ Knebusch}
\newblock{ Weakly Semialgebraic Spaces},
\newblock Lecture Notes in Mathematics 1367, Springer-Verlag, 1989.

\bibitem[{Mcleary01}]{Mcleary}
{\sc J. \ McCleary}
\newblock A User's Guide to Spectral Sequences, Second Edition 
\newblock Cambridge Studies in Advanced Mathematics, 2001.


\bibitem[{Milnor64}]{Milnor}
{\sc J. \ Milnor},
\newblock {\em On the Betti numbers of real varieties},
\newblock Proc. AMS 15, 275-280 (1964).

\bibitem[{Ole\u{i}nik51}]{O}
{\sc O.\ A.\ Ole\u{i}nik},
\newblock {\em Estimates of the {B}etti numbers
of real algebraic hypersurfaces},
\newblock { Mat. Sb. (N.S.)}, 28 (70): 635--640 (Russian) (1951).


\bibitem[{OP49}]{OP}
{\sc O. A.\ Ole\u{i}nik, I. B.\ Petrovskii},
\newblock {\em On the topology of real algebraic surfaces},
\newblock Izv. Akad. Nauk SSSR 13, 389-402 (1949).


\bibitem[{Renegar92}]{R92}
{\sc J.\ Renegar.}
\newblock {\em On the computational complexity and geometry of the
first order theory of the reals},
\newblock Journal  of Symbolic Computation, 13: 255--352 (1992).

\bibitem[{Spanier}]{Spanier}
{\sc E.\ H. Spanier}
\newblock {Algebraic Topology},
\newblock McGraw-Hill Book Company, 1966.

\bibitem[{Rotman}]{Rotman}
{\sc J.\ J. Rotman}
\newblock {An Introduction to Algebraic Topology},
\newblock Springer Verlag, 1988.

\bibitem[{Thom65}]{Thom}
{\sc R.\ Thom},
\newblock {\em Sur l'homologie des vari\'et\'es alg\'ebriques
r\'eelles},
\newblock  {Differential and Combinatorial
Topology},  255--265.
Princeton University Press, Princeton (1965).
\end{thebibliography}
\end{document}